# Flexural Behavior and Design of Prestressed Ultra-High Performance Concrete (UHPC) Beams: Failure Mode and Ductility

Xin TIAN[1], Zhi FANG[2,*], and Yi SHAO[3,*]

## ABSTRACT

Ultra-high performance concrete (UHPC) is well-known for its ultra-high compressive strength and sustained post-cracking tensile ductility, making it an attractive choice for the construction of modern structures. Prestressed UHPC members, however, often fail quickly after crack localization accompanied by reinforcement rupture, which shows limited failure warnings (i.e., relatively small ductility, nearly invisible cracking, and negligible compressive damage). On the other hand, UHPC members can also be designed to allow the gradual strain hardening of reinforcement to compensate for the load loss due to crack localization, and the final failure is attributed to gradual crushing of UHPC prior to reinforcement rupture. Failure after gradual strain hardening is desirable since it brings warning signs through high ductility, visible cracks, and controlled spalling. To achieve a resilient structural design of UHPC members, this study aims to develop design methods to avoid early failure after crack localization and promote failure after gradual strain hardening. This study first establishes a three-dimensional finite-element analysis (FEA) model to simulate the flexural behavior of prestressed UHPC beams, which is validated against seven existing experimental beams. Then, parametric analyses are conducted to assess the influence of five key aspects, including the steel rebar-to-prestressing strand ratio, the post-yield hardening of mild steel rebars, the pretensioned stress of prestressing strands, the tensile behavior of UHPC, and the web width to bottom flange ratio. Additionally, the authors propose a design method to predict the failure mode of prestressed UHPC beams, which can be used to guide the design of reinforcement configuration and to promote ductile structural behavior. Especially, a threshold reinforcing ratio for mild-steel rebar is proposed to design for failure after gradual strain hardening.

[1]Ph.D. Candidate, College of Civil Engineering, Hunan University, Changsha 410082, China. Graduate Research Trainee, Department of Civil Engineering, McGill University, Montreal, Quebec, Canada. (txgone@hnu.edu.cn) ORCID: https://orcid.org/0000-0002-4535-7578.
[2]Professor, College of Civil Engineering, Hunan University, Changsha 410082, China; Key Laboratory for Wind and Bridge Engineering of Hunan Province, Changsha 410082, China. (**Corresponding author:** fangzhi@hnu.edu.cn) ORCID: https://orcid.org/0000-0002-6279-135X.
[3]Assistant Professor, Department of Civil Engineering, McGill University, Montreal, Quebec, Canada. (**Corresponding author:** yi.shao2@mcgill.ca) ORCID: https://orcid.org/0000-0001-9722-9220.

## INTRODUCTION

Ultra-high performance concrete (UHPC) is a class of cementitious composite with ultra-high mechanical strength (a compressive strength of over 120 MPa and a tensile strength greater than 7 MPa (Lagier et al. 2016; ASTM C1856/C1856M-17; Fang et al. 2023)), a high ductile capacity with a tensile localization strain often greater than 0.0025 (Wille et al. 2014), and low-to-negligible permeability (Shi et al. 2015). Among the various innovations, UHPC has been widely applied to form lightweight and high-strength structural elements (Yang et al. 2010; Chen and Graybeal 2012a; Shao and Billington 2022; Shao et al. 2023).

As will be reviewed in the section "Background", many researchers have investigated the flexural behavior of both mild-steel-reinforced and prestressed UHPC beams, which usually experience a brittle failure, namely, reinforcement rupture after crack localization (Shao et al. 2021; Fang et al. 2023). Recently, a new failure mode has been identified through tension-stiffening tests and achieved in both mild-steel-reinforced and prestressed UHPC beams. In this new failure mode, UHPC beams achieve higher load capacity after crack localization and fail by gradual crushing of UHPC after gradual strain hardening of reinforcement. However, this new failure mode has not been fully understood in prestressed UHPC beams, and a design method to predict different failure modes is not established.

Besides experiments, finite-element analysis (FEA) is commonly used to simulate the structural response of UHPC elements, which efficiently reveals the influence of different variables (Lagier et al. 2016; Shafieifar et al. 2018; Shao et al. 2021; Gholipour and Muntasir Billah 2022; Saqif et al. 2022; Shao and Billington 2022; Hu et al. 2023; Peng et al. 2023). Among them, only two-dimensional finite element models with plane-stress elements have been developed to capture the gradual crushing behavior of mild-steel-reinforced UHPC beams with a rectangular cross-section (Shao et al. 2021; Shao and Billington 2022). Three-dimensional finite element models have not been developed for capturing the gradual crushing behavior of UHPC beams with more complex cross-section (e.g., I beam), which is commonly used for prestressed UHPC beams. The recognition of UHPC's crushing behavior as well as a dependable FEA model is vital for the design of prestressed UHPC beams with a more ductile failure mode, namely, failure after gradual strain hardening of steel reinforcement.

To fill the above research gaps, this study aims to establish a robust three-dimensional FEA model to simulate the flexural behavior of prestressed UHPC beams by using the concrete damage plasticity (CDP) model in ABAQUS. First, an FEA model is proposed and verified by seven collected experimental beams in previous studies (El-Helou and Graybeal 2022; Li et al. 2022; Shao and Billington 2022; Fang et al. 2023), highlighting its capacity to capture the load-deflection response as well as the maximum load capacity, crack pattern, and failure mode. Then, parametric analyses are conducted to assess the influence of key parameters on the behavior of prestressed UHPC beams. Finally, the ductility of prestressed UHPC beams is evaluated by displacement ductility factor $\mu$, and a design method is proposed to guide the design of prestressed UHPC beams with a ductile behavior.

## BACKGROUND

Numerous studies have been conducted to address the effect of different parameters (e.g., UHPC strength, steel fiber volume, reinforcing ratio, etc.) on the flexural behavior of reinforced

UHPC (R-UHPC) beams (Yang et al. 2010; Yoo and Yoon 2015; Yoo et al. 2016; Singh et al. 2017; Hasgul et al. 2018; Qiu et al. 2020; Feng et al. 2021; Shao et al. 2021; Shao and Billington 2021, 2022; Qiu et al. 2022). From this extensive experimental database on R-UHPC beams, it is observed that R-UHPC beams have two different failure paths (shown in **Fig. 1**) after crack localization, which show different load-reduction mechanisms and distinct behavior at failure (Shao and Billington 2019). Shao (2019) distinguishes two different failure paths by their different failure mechanisms during loading process: the first failure path is characterized by the loss of the load capacity after crack localization due to fiber-bridging loss, which is accompanied by continuous fiber pull-out and reinforcement rupture at the end; in the second failure path, the beams carry increased external load after crack localization since gradual strain hardening of reinforcement compensates for the loss of fiber-bridging. The second path is recommended because it brings higher deformation ductility, more warning signs before collapse, and full utilization of the compressive property of UHPC.

For an enhanced understanding of this phenomenon, Zhu et al. (2020) tested R-UHPC beams with different steel fiber volumes (i.e., 1.5%, 2.0%, and 3.0%) and reinforcing ratios (1.48%, 2.33%, and 2.83%). They propose that the failure path is strongly dependent on the coupled effect of steel fiber volume and reinforcing ratio and verified one important hypothesis behind the failure path recognition that the crack localization and yielding of mild steel rebars occur simultaneously. Recently, Shao and Billington (2022) and Peng et al. (2023) conducted both experimental and numerical research on the ductility and failure mode of R-UHPC beams. They found that a low reinforcing ratio (below 2%) and UHPC with high tensile strength (above 8 MPa) tend to bring a brittle failure along with reinforcement rupture after crack localization, while a high reinforcing ratio and UHPC with low tensile strength tend to bring a ductile failure along with gradual strain hardening of longitudinal steel reinforcement.

The use of prestressed reinforcement in flexural members, especially in large-scale components,

can improve the crack resistance and flexural capacity of UHPC members. Some scholars have tested the flexural behavior of prestressed UHPC beams, and the observed failure modes were divided into two categories: first, fully prestressed beams and partially prestressed beams with low reinforcing ratio (range from 0.33% to 1.73%) of mild steel rebars usually fail after crack localization with the rupture of steel strands (Graybeal 2008; Yang et al. 2011a; Chen and Graybeal 2012a; Su et al. 2020; El-Helou and Graybeal 2022; Fang et al. 2023); second, the addition of mild steel rebars (e.g., above 2%) leads to a ductile failure mode with UHPC crushing after gradual hardening of reinforcement (Yang et al. 2011b; Xu and Deng 2014; Li et al. 2022). Although failure paths for R-UHPC beams have been studied, research on the failure mechanics for prestressed UHPC beams is limited, especially the gradual crushing behavior of such beams. The coupled effect of steel fiber volume and reinforcing ratio of mild steel rebar and prestressing strand is little understood. The design method to predict the different failure paths for prestressed UHPC beams is not established.

Many studies have examined the FEA models of R-UHPC beams. These models can be divided into two categories: first, two-dimensional FEA models based on the smeared-crack concept or hybrid-rotating/fixed-crack model were established to predict R-UHPC beams that fail after crack localization (Bandelt and Billington 2018), fail after gradual strain hardening (Shao et al. 2021), fail in above two failure paths (Shao and Billington 2022; Peng et al. 2023), or fail in a shear mode (Hung and El-Tawil 2010; Hung et al. 2013); second, three-dimensional FEA models based on the CDP model or smeared-crack concept were established to predict R-UHPC beams that fail after crack localization (Chen and Graybeal 2012a; b; Hung and Li 2013; Mahmud et al. 2013; Shafieifar et al. 2018; Yin et al. 2019) or shear (Bahij et al. 2018; Hussein and Amleh 2018). However, a three-dimensional FEM model that can capture the gradual crushing behavior of UHPC has not been established, which is needed for understanding the behavior of UHPC beams with non-rectangular cross-sections (such as prestressed UHPC I girders).

## MODELING METHODOLOGY

This section presents a numerical model for simulating both mild-steel-reinforced and prestressed UHPC beams. Three-dimensional finite element models were built in ABAQUS, which can capture the nonlinear tensile behavior of UHPC structures (Zhu et al. 2020). **Fig.2** shows the mesh discretization and geometry of the finite-element model for a typical beam with an I-shaped section.

In the actual beams, crack localization commonly initiated at the weakest location in the constant moment region due to initial flaws. To simulate this phenomenon, an initial flaw was introduced in the present simulation by reducing the tensile strength, tensile fracture energy, and compressive fracture energy of three element strips at the midspan by 15% (**Fig. 2**), which is recommended by Shao and Billington (2022).

UHPC was modeled using a three-dimensional eight-node element (C3D8R) with an element size of $10 \times 10 \times 10$-mm. The well-established Concrete Damaged Plasticity (CDP) model was used to simulate the nonlinear behavior of UHPC in tension in ABAQUS (Zhu et al. 2020). The CDP model can handle the transformation of strain from an elastic state to an elastoplastic state and arrive at a fully plastic state in the end. The fundamentals for the plastic behavior of concrete by the CDP model consist of three parts. The first is the tensile behavior of UHPC after cracking, which can be defined by introducing tensile stresses versus inelastic tensile strains. The second is the compressive behavior of UHPC, which can be defined by yielding stresses and inelastic strains. The third is the corresponding tension and compression damage parameters along with the failure surface of concrete elements, determined by five plastic parameters, i.e., flow potential eccentricity ($M$), the ratio of the second stress invariant on the tensile meridian to that on the compressive meridian ($K_C$), the ratio of the initial equiaxial compressive yield stress to the initial uniaxial compressive yield stress ($f_{bo}/f_{co}$), viscosity ($\mu$), and dilation angle ($\phi$). The material models are introduced in the section "MATERIAL CONSTITUTIVE MODELS".

The longitudinal mild steel rebars and prestressed strands were modeled by a 2-node linear 3-D truss element. This element is a long, slender structural member that only transmits axial load. The steel elements were embedded in the UHPC elements, resulting in perfect bonding simulated between the steel reinforcements and UHPC. This assumption simplifies the simulation and considers the superior bonding between steel reinforcement and UHPC (e.g., higher than 30 MPa, (Shao and Ostertag 2022)). The prestressing load was achieved using the prestress field function embedded in ABAQUS.

Loading and support plates were modeled using elastic steel properties with C3D8R elements. The surface-to-surface contact was used for the interface between the loading and support plates and the beam. The hard-contact type was used which allows separation after deformation, as well as tangential behavior with a friction coefficient of 0.30 (Jabbar et al. 2023). The movement and rotation constraints corresponding to the simply-supported boundary conditions are set for support plates.

**Loading and analysis procedure**

Loading and support boundaries were modeled to simulate the conditions in the simply-supported beam under four-point bending. A nonlinear static analysis with an incremental displacement-based loading was used. The step size for each incremental displacement was 0.25 mm.

In the simulation, the fracture of tensile mild steel rebar and prestressing strand occurs when its averaged tensile strain on several integration points exceeds the ultimate tensile strain ($\varepsilon_{sf}$ and $\varepsilon_{pu}$, respectively) (Bandelt and Billington 2018). Crushing of UHPC is assumed when the compressive strain in a single UHPC element reaches the strain corresponding to the stress when it drops to 85% of the peak strength.

# MATERIAL CONSTITUTIVE MODELS

**Constitutive model of UHPC**

The tensile behavior of UHPC was modeled by a tri-linear model (**Fig. 3a**, (Shao et al. 2021; Peng et al. 2023)):

$$\sigma = \begin{cases} E_c \varepsilon & \text{for } 0 \leq \varepsilon \leq \varepsilon_{t,cr} \\ f_{t,cr} + \dfrac{(\varepsilon - \varepsilon_{t,cr})}{(\varepsilon_{t,p} - \varepsilon_{t,cr})} \times (f_{t,p} - f_{cr}) & \text{for } \varepsilon_{t,cr} \leq \varepsilon \leq \varepsilon_{t,p} \\ \dfrac{(\varepsilon_{t,u} - \varepsilon)}{(\varepsilon_{t,u} - \varepsilon_{t,p})} \times f_{t,p} & \text{for } \varepsilon_{t,p} < \varepsilon \leq \varepsilon_{t,u} \end{cases} \quad (1)$$

where $\varepsilon_{t,cr} = f_{t,cr} / E_c$; $\sigma$ = stress, MPa; $\varepsilon$ = strain, mm / mm; $E_c$ = elastic modulus, MPa; $f_{t,cr}$ is cracking tensile strength of UHPC, which is suggested to be $0.95 f_{t,p}$ following Shao et al. (2021); $f_{t,p}$ is peak tensile strength of UHPC; $\varepsilon_{t,cr}$ is UHPC cracking strain; $\varepsilon_{t,p}$ is UHPC localization strain; and $\varepsilon_{t,u}$ is UHPC ultimate tensile strain at which the tensile stress is reduced to zero. The peak tensile strength $f_{t,p}$ and localization strain $\varepsilon_{t,p}$ are the stress and strain where the stress starts to continuously decrease with the increase of the strain, which can be determined from a material test. In the absence of the $\varepsilon_{t,p}$, a minimum localization strain of 0.0025 suggested by El-Helou et al. (2022) was used.

The strain $\varepsilon_{t,u}$ is calculated based on the tensile fracture energy (Shao et al. 2021), $G_f$.

$$\varepsilon_{t,u} = 2 \frac{G_f}{H} \frac{1}{f_{t,p}} + \varepsilon_{t,p} - \frac{f_{t,p}}{E_c} \quad (2)$$

where $G_f$ = tensile fracture energy of UHPC; and $H$ = crack bandwidth. The crack bandwidth depends on the element size and shape, the crack orientation within the element, and the integration scheme. For the adopted C3D8R element, $H$ can be determined as a side length of the three-dimensional elements (Pots 1988), i.e., 10 mm.

The tensile fracture energy, $G_f$, for UHPC can be expressed as a linear function of $v_f l_f / d_f$, if the value of it is not measured or provided in the references. Based on regression by Peng et al.

(2023), the tensile fracture energy can be estimated,

$$G_f = \lambda(11.553\frac{l_f}{d_f}v_f + 5.7859) \tag{3}$$

where $\lambda$ = reduction parameter that accounts for the group fiber orientation, and taken as 0.8 suggested by Peng et al. (2023); $v_f$ = fiber volume fraction; $l_f$ = the length of steel fiber; $d_f$ = the diameter of steel fiber.

Steel fibers can hold the integrity of UHPC matrix in compression, which leads to gradual compression softening and increases its ultimate compressive strain capacity. To simulate the gradual crushing behavior of UHPC, Shao et al. (2021) proposed a fracture-energy-based compression model with a linear descending branch (**Fig. 3b**):

$$\sigma = \begin{cases} E_c\varepsilon & \text{for } \varepsilon_{c,1/2} \leq \varepsilon \leq 0 \\ -\frac{1}{2}f_c' \times \left[1 + k_1 \times \left(\frac{\varepsilon - \varepsilon_{c,1/2}}{\varepsilon_{c,p} - \varepsilon_{c,1/2}}\right) + k_2 \times \left(\frac{\varepsilon - \varepsilon_{c,1/2}}{\varepsilon_{c,p} - \varepsilon_{c,1/2}}\right)^2\right] & \text{for } \varepsilon_{c,p} \leq \varepsilon \leq \varepsilon_{c,1/2} \\ -f_{res} + \frac{\varepsilon_{c,u} - \varepsilon}{\varepsilon_{c,u} - \varepsilon_{c,p}} \times (f_{res} - f_c') & \text{for } \varepsilon_{c,u} \leq \varepsilon \leq \varepsilon_{c,p} \\ -f_{res} \times \frac{51\varepsilon_{c,u} - \varepsilon}{50\varepsilon_{c,u}} & \text{for } \varepsilon \leq \varepsilon_{c,u} \end{cases} \tag{4}$$

where $f_c'$ = compressive strength, MPa; $f_{res}$ = residual compressive strength, typically taken as $0.2f_c'$, MPa; $\varepsilon_{c,1/2}$ = UHPC compressive strain at half of the peak strength, mm / mm; $\varepsilon_{c,p}$ = UHPC compressive strain at the peak strength, mm/mm; and $\varepsilon_{c,u}$ = UHPC compressive strain when decreasing to the residual strength, mm/mm; $k_1 = \frac{\varepsilon_{c,p} - \varepsilon_{c,1/2}}{\varepsilon_{c,1/2}}$; $k_2 = \frac{2\varepsilon_{c,1/2} - \varepsilon_{c,p}}{\varepsilon_{c,1/2}}$.

The UHPC compressive strain at the peak strength calculated by Eq. 5,

$$\varepsilon_{c,p} = -0.724\frac{f_c'}{E_c} - 0.0016 \tag{5}$$

The strain at which the UHPC compressive stress is reduced to $f_{res}$ is calculated by Eq. 6,

$$\varepsilon_{c,u} = -\frac{2G_c}{(f_c' + f_{res})H} + \varepsilon_{c,p} + \frac{f_c' - f_{res}}{E_c} \tag{6}$$

where $G_c$ = compressive fracture energy of UHPC, which is estimated as $\frac{f_c'}{MPa} \times 0.9$ N/mm (Peng et al. 2023) for if $G_c$ is not provided in the literature.

The elastic modulus of UHPC was estimated by Shao et al. (2021):

$$E_c = 3680\sqrt{f_c'} \tag{7}$$

**Damage parameters of UHPC**

To capture the material softening behavior and related damage, the damage parameter is implemented after the peak strength is achieved. The compressive and tensile damage parameters (i.e., $d_c$ and $d_t$) used in the model were developed by other scholars (Birtel and Mark 2006; Zhu et al. 2020):

$$d_k = \frac{(1-\alpha)\varepsilon_{in}E_c}{\sigma_0 + (1-\alpha)\varepsilon_{in}E_c}, \quad (k=t,c) \tag{8}$$

$$\varepsilon_{in} = \varepsilon_k - \frac{\sigma_k}{E_c}, \quad (k=t,c) \tag{9}$$

where $t$ and $c$ mean tension and compression, respectively; $\alpha$ is the scaling coefficient between plastic strain and inelastic strain, which is suggested as 0.35 to 0.70 for compression and 0.50 to 0.95 for tension; $\varepsilon_{in}$ is the inelastic strain for UHPC under tension and compression; $\sigma_0$ is the stress after which the constitutive model entering the plastic stage.

Additionally, the CDP model requires five additional parameters and can be divided into two groups, namely, empirical constants which have little effect on finite element simulation and trial parameters which have a significant role in simulating the flexural behavior of concrete beams in ABAQUS by determining the post-peak performance of concrete (Jabbar et al. 2021; Mahdi 2023).

Empirical constants include flow potential eccentricity ($M$), the ratio of the second stress invariant on the tensile meridian to that on the compressive meridian such that the maximum principal stress is negative ($K_C$), the ratio of the initial equiaxial compressive yield stress to the initial uniaxial compressive yield stress ($f_{bo}/f_{co}$). These parameters were defined as 0.1, 0.6667, and 1.16, respectively. Jabbar et al. (2018) indicated that the higher viscosity ($\mu$) and dilation angle ($\phi$) can lead to a more ductile behavior of UHPC after cracking which makes the simulated beam support higher loads and deform more. Additionally, these two parameters also control the calculation speed and convergence rate. After trial calculations, the viscosity ($\mu$) was defined as 0.0001 in present study. The dilation angle was set as 35°.

**Constitutive model of steel reinforcement**

For mild steel rebars, its stress-strain relationship can be divided into four parts: elastic, yielding, strain hardening, and softening, as shown in **Fig. 4a**. The following Eq. 10 was used for the mild steel rebars (Holzer et al. 1975):

$$f_s = \begin{cases} E_s \varepsilon_s & \text{for } 0 \leq \varepsilon_s \leq \varepsilon_y \\ f_y & \text{for } \varepsilon_y \leq \varepsilon_s \leq \varepsilon_{sh} \\ f_y \left[ 1 + \frac{\varepsilon_s - \varepsilon_{sh}}{\varepsilon_{su} - \varepsilon_{sh}} \left( \frac{f_{su}}{f_y} - 1 \right) \exp\left( 1 - \frac{\varepsilon_s - \varepsilon_{sh}}{\varepsilon_{su} - \varepsilon_{sh}} \right) \right] & \text{for } \varepsilon_{sh} < \varepsilon \leq \varepsilon_{sf} \end{cases} \quad (10)$$

where $f_s$ and $\varepsilon_s$ = stress and strain in steel, respectively; $E_s$ = modulus of elasticity of steel, MPa; $f_y$ and $\varepsilon_y$ = yield stress and strain, respectively; $\varepsilon_{sh}$ = hardening strain of steel; $f_{su}$ and $\varepsilon_{su}$ = maximum stress and corresponding strain of steel, respectively; and $\varepsilon_{sf}$ = rupture strain. The parameters $\varepsilon_{sh}$, $\varepsilon_{su}$, and $\varepsilon_{sf}$ are supposed to be prescribed. In the absence of these parameters, the values proposed by Naaman (2009) were used.

The stress-strain relationship of the prestressing steel strand was that from Mattock (1979) given by the hyperbolic expression valid only in tension, as shown in **Fig. 4b**.

$$f_{ps} = E_{ps} \varepsilon_{ps} \left[ Q + \frac{(1-Q)}{\left( 1 + \left( \frac{E_{ps} \varepsilon_{ps}}{K f_{py}} \right)^N \right)^{1/N}} \right] \leq f_{pu} \quad (11)$$

where $f_{ps}$ and $\varepsilon_{ps}$ are stress and strain in prestressing steel, respectively; $f_{pe}$ is the effective stress of prestressing steel; $E_{ps}$ is modulus of elasticity of prestressing steel; $f_{py}$ is yield stress (typically corresponds to a yield strain equal to 0.01) (ASTM A416/A416M-17), and *K*, *Q*, and *N* are non-dimensional coefficients which may be adjusted to improve fit to experimental data. For Grade 270 strands, the values of *K*, *Q*, and *N* are 1.0852, 0.00772, and 7.39 (Naaman 2009), respectively. To satisfy the minimum specified ASTM Standards, $\varepsilon_{pu}$ is assumed equal to 0.04.

## VALIDATION OF NUMERICAL ANALYSIS

The proposed numerical model was validated by comparing the numerical results to experimental results obtained from three fully or partial prestressed UHPC beams (El-Helou and Graybeal 2022; Li et al. 2022; Fang et al. 2023) and four R-UHPC flexural members (Shao and Billington 2022). **Appendix A** summarizes details of the geometry and material properties of the specimens considered. All specimens were reported to fail in flexure. **Table 1** lists the material properties of the above UHPC beams. The load-deflection curves obtained from the numerical model of the seven beams are presented in **Fig. 5** along with the experimentally-measured load-deflection curves.

As can be seen, the load-deflection response obtained from the proposed numerical model agrees well with the experimental results including the initial phases, crack localization points, and two failure paths. The final load decline before fracture was not well captured, which may be attributed to the deviation between the assumed material model and actual behavior (i.e., the necking of steel reinforcement which may bring abrupt reduction of the curve).

**Fig. 6** also compares the predicted and experimental crack patterns. Overall, the numerical results accurately capture the cracking behavior. In this figure, the Damage T and Damage C represent the tensile damage parameter $d_t$ and the compressive damage parameter $d_c$, respectively. As seen from **Fig. 6 (a), (b), and (f)**, a localized crack ($d_t \geq 0.90$) appeared near the loading point with the concentration of plastic strain of steel strands around the localized crack for specimens failed due to fracture of tensile steel reinforcement after crack localization. For **Fig. 6 (c), (d), and (e)**, fine and dense cracks formed in the constant moment region of simulated specimens with a longer plastic hinge length of steel plastic strain. Finally, specimens failed after the crushing of UHPC ($d_c \geq 0.60$).

For **Fig. 6 (g)**, a localized crack appeared first but steel strain hardening compensated for the initial loss of fiber-bridging introduced by crack localization. Then, the load-deflection curve remained plateau and the load reduction was ultimately triggered by fiber-bridging loss rather than by crushing of UHPC (Shao and Billington 2022).

## PARAMETRIC ANALYSIS

After the verification of proposed numerical model, a systematic parametric analysis was conducted on prestressed UHPC beams with a prototype of Fang et al.'s experiment (2023) (**Fig. 6a**). Numerical test matrix was created based on the key aspects that may affect the failure modes of the beam including the steel rebar-to-prestressing strand ratio, the tensile behavior of steel reinforcement, the pretensioned stress of strand, the tensile behavior of UHPC (i.e., tensile strength and localization strain), and the web width to bottom flange ratio. **Appendix B** lists all the numerical specimens with their detailed modeling information.

Each specimen was given a name which is divided into three parts. The letter "M" and the number behind "M" denotes the amount of mild steel rebar, and the number behind "@" means the diameter of mild steel rebar. The letter "P" and the number behind "P" denotes the amount of prestressing strand, and the number behind "@" means the diameter of the strand. The last part means the type of parameter analysis. "Re" means reference model. "SH" means post-yield strain hardening strength of mild steel rebars. "$f_{con}$" means the pretensioned stress of prestressing strands. "$f_t$" means the tensile strength of UHPC. "loc" means the tensile localization strain of UHPC. "$b_w$" means the web width of girder. For example, M6@22.2-P6@15.2-SH250 denotes a girder reinforced by six mild steel rebars with a diameter of 22.2 mm and a post-yield strain hardening strength of 250 MPa, and six prestressing strands with a diameter of 15.2 mm. In addition, a reference specimen has a post-yield hardening strength of 150 MPa for mild steel, a pretensioned stress of 1395 MPa for prestressing strand, a tensile strength of 9 MPa for UHPC, a localization strain of 0.25% for UHPC, and a web width of 250 mm.

**Effect of steel rebar-to-prestressing strand ratio**

This section explores the effect of the steel rebar-to-prestressing strand ratio, i.e., different combinations of mild and prestressed steel for the same flexural load capacity. An estimation method is first introduced to guide the design of different specimens with similar load capacity. **Fig. 7** shows the strain and stress distribution across the depth of a typical prestressed UHPC section at crack localization. In the figure, $d_s$ and $d_p$ are the distance between the centroid of mild steel rebars and prestressing strands to the top of the beam, respectively. $\varepsilon_{p,e}$ is the effective initial strain of prestressing strands. $\varepsilon_{s,loc}$ and $\varepsilon_{p,loc}$ are the strains of mild steel rebar and prestressing strand at UHPC localization point, respectively. $\varepsilon_c$ and $f_c$ are the strain and stress at the extreme compression fiber of the section, respectively. $\varepsilon_{t,p}$ and $f_{t,p}$ are the localization strain and peak tensile strength of UHPC, respectively.

The proportion of prestressing strands and mild steel rebars was achieved by the equivalence of load capacity at the crack localization state. The crack localization load capacity is dependent on the total tension force carried by UHPC, mild steel rebar, and prestressing strand. Then the section area of mild steel rebar and prestressing strand can be estimated through the equivalence of tension force.

$$T = T_{uc} + T_{M,loc} + T_{P,loc} \tag{12}$$

$$T = \gamma bh f_t + A_s f_{s,loc} + A_p f_{p,loc} \tag{13}$$

where $T_{uc}$, $T_{M,loc}$, and $T_{P,loc}$ are the total tension force of UHPC, mild steel rebar, and prestressing strand, respectively; $f_t$ = UHPC tensile strength; $\gamma$ = empirical factor representing the portion of the cross-sectional area that has fiber-bridging, where $\gamma$ is calibrated to be 0.54 based on 35 reinforced HPFRCC flexural test results from literature (Shao and Billington 2019); $h$ and $b_w$ = beam height and web width; $f_{s,loc}$ and $f_{p,loc}$ = the stress of mild steel rebar and prestressing strand, respectively, when the UHPC reaches crack localization state. $f_{s,loc}$ is assumed to be the yielding stress of mild steel through experimental and numerical investigations (Shao and Billington 2022; Zhang et al. 2023). $f_{p,loc}$ is calculated through plane section assumption as shown in **Fig. 7**; $A_s$ and $A_p$ = the section area of mild steel and prestressing strand, respectively.

Three numerical models were established for reference, namely M2@16-P10@15.2-Re, M6@22.2-P6@15.2-Re, and M10@22.7-P2@15.2-Re. It is worth noting that the simulated rebar diameter is not a common diameter used because this paper focuses on finite element calculation to clarify the influence of varied factors on prestressed UHPC beams. Therefore, the simulated rebar diameter was calculated from the above Eqs. 12 and 13.

**Fig. 8** shows the load-deflection curves and tensile damage contours of prestressed UHPC beams with different arrangements of steel reinforcements. In this figure, the cross-section diagram shows that the different steel rebar-to-prestressing strand ratio is achieved by switching the mild-steel rebar and prestressing strand with the same height of gravity position to keep the location of the reinforcements unchanged. These three numerical beams were designed with similar flexural capacity and achieved a range of drift capacities from 2.87% to 3.26% along with two failure paths. A larger mild steel rebar-to-prestressing strand reinforcing ratio (i.e., $A_s / A_p$) increases the drift capacity of beams. For specimens that failed after crack localization (sudden drop of applied load after point ■), increasing the amount of longitudinal mild steel rebars either changed the failure path to failure after gradual strain hardening (compare M2@16-P10@15.2-Re to M10@22.7-P2@15.2-Re) or slowed the load drop (compare M2@16-P10@15.2-Re to M6@22.2-P6@15.2-Re). The tensile damage contour shows that a larger amount of longitudinal mild steel rebars resulted in larger area of damaged elements (more localized cracks) and softened materials, indicating a larger region of plasticity.

**Effect of post-yield hardening strength of mild steel rebar**

Three values of post-yield hardening strength ($\Delta f_u$ = 50, 150, and 250 MPa) were considered. These simulated beams both had a yield strength of 485 MPa. **Fig. 9** shows the load-deflection response of prestressed UHPC beams with different post-yield hardening strengths of steel. The following conclusions can be found:

A higher post-yield hardening of mild steel rebar is found to increase the drift capacity of the beams by either changing the failure path to failure after gradual strain hardening (compare M6@22.2-P6@15.2-Re to M6@22.2-P6@15.2-SH250) or slowed the load drop from crack localization (compare M2@16-P10@15.2-Re to M2@16-P10@15.2-SH50). When the post-yield hardening strength provided by mild steel rebar is larger than the fiber-bridging capacity provided by steel fibers, the failure mode will change to failure after gradual strain hardening. The reason behind that is the total tension force transferred by the hardened mild steel rebar and prestressing strand is larger than the total tension force carried by the yielded mild steel rebar, prestressing strand, and fiber at the crack localization. This tension force increase leads to a higher load capacity after the initial crack localization , with more localized cracks being able to form. This phenomenon was also numerically found in the tensile damage contour where specimens, e.g., M6@22.2-P6@15.2-SH250 and M10@22.7-P2@15.2-SH250, that failed after gradual strain hardening had more localized cracks whereas the specimens, e.g., M2@22.2-P10@15.2-SH50 and M6@22.7-P6@15.2-SH50, failed after crack localization only had one localized crack. The occurrence of several localized cracks allowed the steel to harden across a longer plasticity region, as indicated by the presence of more red-colored elements and increased the specimen deformation capacity, which was also observed in tension-stiffening tests and flexural tests of mild-steel-reinforced UHPC members (Shao et al. 2021; Zhang et al. 2023).

**Effect of pretensioned stress of prestressing strand**

Two values of pretensioned stress of strand ($f_{con}$ = 930 and 1395 MPa) were considered. **Fig. 10** shows the load-deflection response of prestressed UHPC beams with different pretensioned stress of strand. A lower pretensioned stress of strand is observed to increase the drift capacity of the beams by changing the failure path to failure after gradual strain hardening (compare M6@22.2-P6@15.2-Re to M6@22.2-P6@15.2-fcon930 and M10@22.2-P2@15.2-Re to M10@22.2-P2@15.2-fcon930). This is reasonable since the dropping in pretensioned stress of strand results in lower tensile stress at the crack localization point, i.e., a larger amount of tensile stress can be utilized after crack localization to compensate for the load loss after crack localization. From M2@16-P10@15.2-Re to M2@16-P10@15.2-fcon930, the combined hardening strength is still not able to compensate for the load loss after crack localization, which will be demonstrated through the proposed design method in the section "THRESHOLD REINFORCING RATIO FOR DUCTILE FAILURE". Thus, both M2@16-P10@15.2-Re and M2@16-P10@15.2-fcon930 failed after crack localization.

**Effect of tensile properties of UHPC**

Three values of tensile strength of UHPC ($f_{t,p}$ = 6, 9, and 12 MPa) were considered. These simulated beams both had a localization strain of 0.25% for UHPC. **Fig. 11** shows the load-deflection response of prestressed UHPC beams with different tensile strengths of UHPC. Comparing specimens with similar reinforcement arrangements and different tensile strengths from 6 MPa to 12 MPa is similar to the previous comparison of specimens with similar reinforcement arrangements and different post-yield hardening strengths Increasing the UHPC tensile strength made the tension force transferred by fiber-bridging higher than the post-yield-hardening strength provided by mild steel rebar. As a result, the following outcomes were observed (1) resulted in one single localized crack, (2) allowed the steel to harden across a shorter plasticity region, and (3) decreased the drift capacity. Additionally, the increased tensile strength led to an increased flexural capacity.

Three values of localization strain of UHPC ($\varepsilon_{t,p}$ = 0.1%, 0.25%, and 0.6%) were considered. These simulated beams both had a tensile strength of 9 MPa. **Fig. 12** shows the load-deflection response of prestressed UHPC beams with different localization strain of UHPC. The influence of UHPC localization strain is most impactful to the specimens that fail after crack localization since this failure mode was dominated by UHPC tension behavior. For specimens with low mild steel rebar-to-prestressing strand reinforcing ratios, a larger UHPC localization strain delayed the crack localization (compare M2@16-P10@15.2-loc0.1% to M2@16-P10@15.2-loc0.6%) and the associated load drop (compare M6@16-P6@15.2-loc0.1% to M6@16-P6@15.2-loc0.6%), increasing the drift capacity or the formation of a second localized crack (compare M6@16-P6@15.2-loc0.1% to M6@16-P6@15.2-loc0.6%).

For specimens reinforced with ten mild steel rebars and two prestressing strands, two phenomena were observed from the comparison between M10@16-P2@15.2-loc0.1% to M10@16-P2@15.2-Re and M10@16-P2@15.2-loc0.6% to M10@16-P2@15.2-Re. First, the decreased localization strain led to a longer post-peak plateau. Since the reduction in localization strain lowered the reinforcement tension force/stress at crack localization point, a larger amount of tension force is available to compensate for the load loss after crack localization, resulting in a longer plasticity region and higher drift capacity. Second, the increased localization strain led to a higher load bearing capacity and a pronounced descending portion of the load-deflection curve. This is because the increased localization strain will slightly increase the load capacity but in the later loading period, the remaining post-localization hardening capacity of reinforcement is decreased. Therefore, the hardening of steel cannot compensate for the reduction of UHPC tension, which leads to an abrupt load drop.

It should be noted that, while UHPC shows a much higher tensile ductility than conventional concrete, its tensile strain capacity (e.g., a localization strain of 0.1%-0.6%) is still much lower than the tensile strain capacity of longitudinal reinforcement (e.g., 7%-9% for prestressing strands). Therefore, forming multiple localized cracks and a longer plasticity region allows better utilization of the plasticity of longitudinal reinforcement, which contributes to higher structural ductility in prestressed UHPC beams.

**Effect of web width to bottom flange ratio**

Two values of web width ($b_w$ = 130 and 250 mm) were considered. **Fig. 13** shows the load-deflection response of prestressed UHPC beams with different web widths. A higher web width to bottom flange ratio decreased the drift capacity by changing the failure path to failure after crack localization. This is expected since the larger beam width made a larger contribution to fiber-bridging in tension area, in turn, the post-yield hardening strength cannot compensate for the loss of tensile force provided by UHPC.

## EVALUATION OF DUCTILITY

Ductility can reflect the ability of a structural member to bear inelastic deformation prior to collapse, without significant loss in load resistance (Peng et al. 2023). For structural members, ductility is defined as the ratio of deformations at the ultimate state to those at the first yielding of tensile reinforcements. It is worth noting that prestressing strands do not have a clear yielding platform, a specified yielding strength is set (the minimum yield strength shall be 0.90 $f_{p,u}$) when the yielding strain reaches 0.01 as per general design code (ACI 318–19; ASTM A416/A416M-17). Then, the yielding of prestressing strands is defined when the averaged tensile strain on several integration points exceeds 0.01. The displacement ductility factor $\mu$ is expressed as follows,

$$\mu = \frac{\Delta_u}{\Delta_y} \tag{14}$$

where $\Delta_u$ and $\Delta_y$ = the midspan deflection at the ultimate state and the yielding of steel reinforcement, respectively.

To assess the impact of material properties and reinforcement configuration on the ductility of prestressed UHPC beams, a six-factor two-level full factorial design was conducted using Minitab Version 21.1 (Pokhrel and Bandelt 2019). This design is commonly used to determine the influence of different variables in experiments involving multiple variables. Further, it can demonstrate whether the combined effect of two or more variables is significant in an experiment or not.

The present study considers six factors (steel rebar-to-prestressing strand ratio, post-yield hardening strength of mild steel rebar, pretensioned stress of strand, UHPC tensile strength, UHPC localization strain, and web width to bottom flange ratio). Four factors, including the steel rebar-to-prestressing strand ratio, post-yield hardening strength of mild steel rebar, UHPC tensile strength, and UHPC localization strain, contained three levels, while the other two factors, e.g., pretensioned stress of strand and web width to bottom flange ratio, contained two levels. The statistical analysis, performed at a 95% confidence level using a two-sided confidence interval, was based on the full factorial design feature in Minitab.

The output of the standardized effects in the form of a Pareto chart (**Fig. 14a**) indicates that the steel rebar-to-prestressing strand ratio, post-yield hardening strength of mild steel rebar, pretensioned stress of strand, UHPC tensile strength, and web width to bottom flange ratio cross the reference line. Therefore, these five factors have a significant effect on the ductility factor. **Fig. 14a** shows the degree to which a factor is statistically significant based on the absolute value of the standardized effect (X-axis) of the factors. It can be observed that the steel rebar-to-prestressing strand ratio had the strongest effect, followed by UHPC tensile strength, and then post-yield hardening strength of mild steel rebar. The effects due to UHPC localization strain (or combination of localization strain & steel rebar-to-prestressing strand ratio, etc.) were not significant as the standardized effect was less than the reference value ($< 2.13$). The normal probability plot of the standardized effects shows whether the effect increases or decreases the response (i.e., ductility), as shown in **Fig. 14b**.

The normal plot shows positive standardized effects on the right side and negative standardized effects on the left of the red solid line. In other words, positive values signify that ductility is positively correlated to the factor, while negative values mean that ductility is negatively correlated. Specifically, the effect of pretensioned stress of strand, UHPC tensile strength, and web width to bottom flange ratio (red cube) is negative meaning that the ductility decreases with increasing these three factors. In contrast, ductility increases with beams arranged with more mild steel rebars and increasing post-yield hardening strength of mild steel rebar (black cube).

The factors with insignificant effects are positioned near zero or close to the red line, while the factors shown away from the red line have significant effects on the ductility. The extent to which these factors affect the ductility can be determined by observing their distance from the red line. For example, the steel rebar-to-prestressing strand ratio (Point F) causes a sharp increase in ductility and is therefore furthest away from the red line (right). UHPC localization strain (Point D) has the smallest effect and is therefore near zero. The increasing tensile strength (Point C) leads to the largest decrease in ductility for prestressed UHPC beams.

## THRESHOLD REINFORCING RATIO FOR DUCTILE FAILURE

Shao and Billington (2019) proposed a ratio ω to evaluate the drift capacity (defined as the ratio of the midspan deflection to the shear-span length) of mild-steel-reinforced UHPC beams based on three parameters, namely, reinforcing ratio, post-yield hardening strength of mild steel rebar, and UHPC tensile strength.

$$\omega = \frac{(f_{su} - f_{s,loc})\rho_s}{\gamma f_t} \frac{d}{h} \qquad (15)$$

where $f_{su}$ = ultimate strength of the mild steel rebar; $d$ and $h$ = beam effective depth and height; $\rho_s$ = the reinforcing ratio of mild steel rebar.

The strain-hardening portion provided by prestressing strands is added to the Eq. 15.

$$\omega = \frac{\left(f_{su} - f_{s,loc}\right)\rho_s + \left(0.94 f_{pu} - f_{p,loc}\right)\rho_p}{\gamma f_t} \frac{d}{h} \quad (16)$$

where $f_{pu}$ = ultimate strength of the prestressing strand. 0.94 is a reduction factor due to the brittle fracture of steel strand according to ACI 318-19 (2019); $\rho_p$ = the reinforcing ratio of prestressing strand.

**Fig. 15** lists the experimentally and numerically observed displacement ductility ratio $\mu$ with $\omega$. When $\omega < 1$ (the area at the left side of the blue dashed line in the figure), the beams mostly fail after crack localization or fail in a transition state (i.e., load initially increases after crack localization but quickly reduces due to fiber-bridging capacity loss) with lower displacement ductility because the steel hardening capacity cannot compensate for the loss of fiber-bridging capacity. When $\omega \geq 1$ (the area at the right side of the blue dashed line in the figure), the specimens mostly fail after gradual strain hardening with higher displacement ductility because the steel hardening capacity is stronger than the fiber-bridging capacity.

To avoid failure after crack localization, the designer can design for ω larger than 1, which is also supported by **Fig. 15**. Additionally, several simplifications could be applied to develop an easy-to-use equation for practice engineers:

- For mild steel, the localization point stress $f_{s,loc}$ could be taken as $f_{sy}$. This is because the localization strain is typically around the yielding plateau of mild steel. Further, In ACI 318-19 (2019), the tensile strength of the steel $f_{su}$ is taken as the minimum value of 1.25 $f_{sy}$. Therefore, $f_{su} - f_{s,loc}$ can be conservatively simplified as $0.25 f_{sy}$.

- For prestressing strand, $f_{pu}$ can be taken as the design tensile strength of prestressing strand. $f_{p,loc}$ can be estimated using Eq. 11 while the strand strain at localization point could be estimated as the addition of effective prestressing strain and UHPC localization strain.

- ACI 318-19 (2019) also points out that the effective depth *d* need not be taken less than 0.8*h*. *d* is assumed to be 0.8*h* for a conservative design.

With $\omega \geq 1$, Eq. 16 can be rearranged and simplified into:

$$\rho_s \geq \frac{0.675 f_t - (0.94 f_{pu} - f_{p,loc})\rho_p}{0.25 f_{sy}} \quad (17)$$

When designing a prestressed beam, the prestressing reinforcement is typically designed to avoid tensile cracking under service loads. Therefore, the designers can use the developed mild steel threshold equation (i.e., Eq. 17) to help target a ductile behavior with the selected UHPC material and prestressing level.

## FLEXURAL STRENGTH PREDICTION

### Calculation of the ultimate flexural capacity, $M_u$

When the ultimate tensile strength of steel strand $f_{pu}$ is reached, the UHPC tensile resistance can be ignored due to the wide-opening localized crack (Fang et al. 2023). As shown in **Fig. 16**, the ultimate flexural capacity, $M_u$, can be calculated by following equation,

$$M_u = A_p f_{pu}(d_p - c_u) + A_s f_{su}(d_s - c_u) + \sum f_{ci} b_i (x_i - c_u) \Delta h \quad (18)$$

$$\sum f_{ci} b_i \Delta h = A_p f_{pu} + A_s f_{su} \quad (19)$$

$$\frac{c_u}{d_p - c_u} = \frac{\varepsilon_{cf}}{\varepsilon_{pu} - \varepsilon_{pe}} \quad (20)$$

where $c_u$ is depth of the neutral axis at ultimate state; $d_s$ and $d_p$ are the distance between mild steel rebar and prestressing strand to the top of the beam; and $f_{su}$ is ultimate stress in mild steel rebar, and can be taken as 1.25 $f_{sy}$ specified in ACI 318-19 (2019); $f_{ci}$ is UHPC compressive stress at the centroid of the *i*-th layer; $b_i$ is width of the *i*-th layer; $\Delta h$ (= *h* / *n*) is thickness of the *i*-th layer; $x_i$ is distance from the mid-depth of *i*-th segment to the extreme compression fiber.

### Calculation of the crack localization moment, $M_{loc}$

**Fig. 17** shows the strain and stress distributions across the depth of a prestressed UHPC section at crack localization. In this case, the strain in steel strands reaches their yield strain, $\varepsilon_{py}$, while the UHPC in the compression zone of the section is usually in the elastic stage. The nonlinear tensile stress distribution in UHPC can be approximated with an equivalent rectangular stress block using two parameters, $\alpha$ and $\beta$. Then, the moment at crack localization can be determined by the following equations:

$$M_{loc} = \alpha f_{t,p}(b-b_w)\left[h_f - (1-\beta)(h-c)\right]\left(h-c-\frac{h_f - (1-\beta)(h-c)}{2}\right) + \alpha f_{t,p} b_w \beta(h-c)\frac{\beta(h-c)}{2}$$
$$+ A_p f_{py}(d_p - c) + A_s f_{sy}(d_s - c) + \sum f_{ci} b_i (x_i - c)\Delta h$$

(21)

$$\sum f_{ci} b_i \Delta h = \alpha f_{t,p}(b-b_w)\left[\beta(h-c) - (h-c-h_f)\right] + \alpha f_{t,p} b_w \beta(h-c) + A_p f_{py} + A_s f_{sy} \quad (22)$$

$$\frac{c}{d_p - c} = \frac{\varepsilon_{cf}}{\varepsilon_{py} - \varepsilon_{pe}} \quad (23)$$

where $\varepsilon_{cf}$ is UHPC strain at the extreme compression fiber of the section; $b_w$ is width of web (width of a rectangular cross-section); $b$ is width of tension edge of member; and $h_f$ is depth of the bottom flange.

In previous studies (Peng et al. 2022; Fang et al. 2023), the authors suggested the parameter $\alpha$ might be taken as 0.9. The parameter $\beta$ is found to depend on the section height, and can be determined by,

$$\beta = 1.12 - \varsigma \frac{h}{l_f} \leq 1.0 \quad (24)$$

where $\varsigma$ is a factor that depends on the yield strain of steel and prestressing level. For Grade 270 strands with $f_{pe} \geq 0.5 f_{pu}$, $\varsigma$ can be taken as 0.004.

**Validation of the flexural strength prediction method**

The proposed method was already evaluated by the parallel experimental database of prestressed UHPC beams established by the authors (Fang et al. 2023). For further promotion, the numerical results obtained from the current study were also used for the evaluation of this method. Besides, the predictions by the FHWA (Graybeal and El-Helou 2023) were also compared with the proposed model to further evaluate its performance.

**Table 2** compares the numerical flexural strength prediction for prestressed UHPC beams at crack localization and ultimate state to the predicted strength. It was found that the proposed method can well predict the flexural strength of prestressed UHPC beams. At the crack localization state, the mean and standard deviation of the predicted-to-numerical measured strength ratio are 0.97 and 0.10, respectively. At the ultimate state, the mean and standard deviation of the predicted-to-numerical measured strength ratio are 1.02 and 0.11, respectively. The FHWA method, however, tends to give a conservative prediction for the ultimate flexural strength of prestressed UHPC beams with the mean and standard deviation of the predicted-to-numerical measured strength ratio are 0.86 and 0.16, respectively. The conservative prediction is likely due to the way they treat compression and tension stress blocks which is a main difference between calculation methods based on sectional analysis: (1) FHWA assumed that the ultimate capacity of beam coincides with the compressive failure of UHPC when the strain reaches the ultimate compressive strain. Differently, the present model assumed that the ultimate state is accompanied with the rupture of longitudinal reinforcement; (2) they introduce a reduction factor 0.85 on compressive strength to reflect the linearity limit of the compressive behavior of UHPC; (3) they introduce a reduction factor which is not to exceed 0.85 on tensile strength to account for variability in the material behaviors and address the fiber orientation effects which may cause reduction in tensile property (Islam et al. 2022; Zhan et al. 2023).

## CONCLUSIONS

This study investigates the failure mode and ductility of prestressed UHPC beams. A three-dimensional FEA model based on ABAQUS is used to simulate the flexural behavior of prestressed UHPC beams, which is validated against seven existing experimental beams. Then, parametric analyses are conducted on 27 numerical simulations to assess the influence of the steel rebar-to-prestressing strand ratio, the post-yield hardening of mild steel rebar, the pretensioned stress of prestressing strand, the tensile behavior of UHPC, and the web width to bottom flange ratio. Additionally, the authors propose a design method of reinforcement configuration based on the ratio between the steel post-yield hardening capacity and fiber-bridging capacity to ensure a ductile failure, as well as a threshold reinforcing ratio for mild-steel rebar based on the calculation of the ratio $\omega$. The following conclusions are reached:

- Two failure paths and their transition state are observed in numerical simulations: (1) beams that fail after crack localization lose load capacity after crack localization accompanied by reinforcement rupture, which shows limited failure warnings (i.e., relatively small ductility, nearly invisible cracking, and negligible compressive damage); and (2) beams that fail after gradual strain hardening with increasing load capacity until UHPC crushes and show large deformation capacity. A transition state exists where after crack localization, beams may carry additional load by hardening of steel but soon lose load carrying capacity due to fiber-bridging decline instead of reaching the crushing of UHPC.

- Through numerical studies, it was observed that configuring more longitudinal mild steel rebars, lowering UHPC tensile strength, increasing reinforcing steel post-yield hardening strength, and lowering the web width to bottom flange ratio are effective ways to improve the ductility and to change the failure path of flexural elements from crack localization to failure after gradual strain hardening. Changing the pretensioned stress of strand and UHPC localization strain (from 0.1% to 0.6%) was most impactful to the specimens that failed after crack localization by changing their failure mode and only influenced the post-peak portion of load-deflection curves for specimens that failed after gradual strain hardening.

- The ratio between the steel post-yield hardening capacity and fiber-bridging capacity ($\omega$) is an accurate way to guide the design of reinforcement configuration to ensure a ductile failure. When $\omega < 1$, the beams mostly fail after crack localization or fail in a transition state because the steel hardening capacity cannot compensate for the loss of fiber-bridging capacity. When $\omega \geq 1$, the specimens mostly fail after gradual strain hardening because the steel hardening capacity is stronger than the fiber-bridging capacity.

- A threshold reinforcing ratio for mild-steel rebar is proposed based on the ratio $\omega$ to help target a ductile behavior with desired UHPC material and prestressing level.

- Comparisons between the predicted and numerical results show that the proposed method provides accurate predictions on the flexural strength of prestressed UHPC girders. At the crack localization state, the mean and standard deviation of the predicted-to-numerical measured strength ratio are 0.97 and 0.10, respectively.

## ACKNOWLEDGEMENTS

The authors gratefully appreciate the financial support from Canadian Precast/Prestressed Concrete Institute (CPCI), NSERC ALLRP 591005-2023, and MITACS IT39607 through the project "Enhancing the Structural Ductility of Prestressed UHPC Girder". The authors gratefully acknowledge the financial support provided by the National Natural Science Foundation of China (Project No. 51938012). The first author also appreciates the support for studying abroad by the China Scholarship Council program (Project No. 202306130103).

## DATA AVAILABILITY STATEMENT

Some or all data, models, or code that support the findings of this study are available from the corresponding author upon reasonable request.

**Table 1**. Mechanical properties of UHPC.

| Material properties | Notation | Unit | Fang et al. | El-Helou and Graybeal | Li et al. | Shao-vf1.0 | Shao-vf2.0 |
|---|---|---|---|---|---|---|---|
| Elastic modulus | $E_c$ | GPa | 53.8 | 43.4 | 42.2 | 53.5 | 53.5 |
| Compressive strength | $f_c'$ | MPa | 140.6 | 173.0 | 128.0 | 185.3 | 185.8 |
| Compressive fracture energy | $G_c$ | MPa-mm | 126.5 | 155.7 | 115.2 | 167.2 | 166.8 |
| Tensile strength | $f_{t,p}$ | MPa | 9.2 | 10.4 | 7.1 | 7.0 | 10.5 |
| Tensile localization strain | $\varepsilon_{t,p}$ | % | 0.25 | 0.50 | 0.25 | 0.02 | 0.20 |
| Tensile fracture energy | $G_f$ | MPa-mm | 16.6 | 16.6 | 20.65 | 10.4 | 11.2 |
| Poisson's ratio | $\upsilon$ | mm/mm | 0.18 | 0.18 | 0.18 | 0.18 | 0.18 |

**Table 2.** Evaluation of the flexural strength prediction method.

| Name | $M_{loc,Cal}$ | $M_{loc,Num}$ | | $M_{u,Cal}$ | $M_{u,Num}$ | | $M_{u,FHWA}$ | |
|---|---|---|---|---|---|---|---|---|
| M2@16-P10@15.2-Re | 2363 | 2399 | 0.99 | 2231 | 2281 | 0.98 | 1940 | 0.85 |
| M6@22.2-P6@15.2-Re | 2271 | 2386 | 0.95 | 2286 | 2293 | 1.00 | 1850 | 0.81 |
| M10@22.7-P2@15.2-Re | 2111 | 2384 | 0.89 | 2252 | 2390 | 0.94 | 1699 | 0.71 |
| M2@16-P10@15.2-SH50 | 2363 | 2436 | 0.97 | 2203 | 2198 | 1.00 | 1940 | 0.88 |
| M6@22.2-P6@15.2-SH50 | 2271 | 2424 | 0.94 | 2123 | 2210 | 0.96 | 1850 | 0.84 |
| M10@22.7-P2@15.2-SH50 | 2111 | 2504 | 0.84 | 1968 | 2142 | 0.92 | 1699 | 0.79 |
| M2@16-P10@15.2-SH250 | 2363 | 2385 | 0.99 | 2259 | 2308 | 0.98 | 1940 | 0.84 |
| M6@22.2-P6@15.2-SH250 | 2271 | 2237 | 1.02 | 2449 | 2539 | 0.96 | 1850 | 0.73 |
| M10@22.7-P2@15.2-SH250 | 2111 | 2385 | 0.89 | 2534 | 2620 | 0.97 | 1699 | 0.65 |
| M2@16-P10@15.2-fcon930 | 2363 | 2500 | 0.95 | 2232 | 2085 | 1.07 | 1940 | 0.93 |
| M6@22.2-P6@15.2-fcon930 | 2271 | 2290 | 0.99 | 2287 | 2460 | 0.93 | 1850 | 0.75 |
| M10@22.7-P2@15.2-fcon930 | 2111 | 2318 | 0.91 | 2253 | 2268 | 0.99 | 1699 | 0.75 |
| M2@16-P10@15.2-ft6 | 2257 | 1846 | 1.22 | 2231 | 1677 | 1.33 | 1922 | 1.15 |
| M6@22.2-P6@15.2-ft6 | 2165 | 1921 | 1.13 | 2286 | 2434 | 0.94 | 1832 | 0.75 |
| M10@22.7-P2@15.2-ft6 | 2005 | 2167 | 0.93 | 2252 | 2402 | 0.94 | 1683 | 0.70 |
| M2@16-P10@15.2-ft12 | 2469 | 2258 | 1.09 | 2231 | 1777 | 1.26 | 1961 | 1.10 |
| M6@22.2-P6@15.2-ft12 | 2377 | 2526 | 0.94 | 2286 | 2056 | 1.11 | 1870 | 0.91 |
| M10@22.7-P2@15.2-ft12 | 2218 | 2611 | 0.85 | 2252 | 2188 | 1.03 | 1717 | 0.78 |
| M2@16-P10@15.2-loc0.1% | 2363 | 2171 | 1.09 | 2231 | 1865 | 1.20 | 1940 | 1.04 |
| M6@22.2-P6@15.2-loc0.1% | 2271 | 2428 | 0.94 | 2286 | 2116 | 1.08 | 1850 | 0.87 |
| M10@22.7-P2@15.2-loc0.1% | 2111 | 2352 | 0.90 | 2252 | 2594 | 0.87 | 1699 | 0.66 |
| M2@16-P10@15.2-loc0.6% | 2363 | 2487 | 0.95 | 2231 | 2448 | 0.91 | 1940 | 0.79 |
| M6@22.2P6@15.2-loc0.6% | 2271 | 2552 | 0.89 | 2286 | 2354 | 0.97 | 1850 | 0.79 |
| M10@22.7P2@15.2-loc0.6% | 2111 | 2464 | 0.86 | 2252 | 2395 | 0.94 | 1699 | 0.71 |
| M2@16-P10@15.2-b250 | 2995 | 2600 | 1.15 | 2231 | 2045 | 1.09 | 2515 | 1.23 |
| M6@22.2P6@15.2-b250 | 2901 | 2663 | 1.09 | 2286 | 2149 | 1.06 | 2415 | 1.12 |
| M10@22.7P2@15.2-b250 | 2734 | 2840 | 0.96 | 2252 | 2096 | 1.07 | 2213 | 1.06 |
| | | Ave. | 0.97 | | Ave. | 1.02 | Ave. | 0.86 |
| | | Std. | 0.10 | | Std. | 0.11 | Std. | 0.16 |

**Appendix A.** Details of geometry and material properties of specimens.

| Investigators | Name | Type | Cross Section | | | UHPC Property | | | | Longitudinal Reinforcement | | | | | | Failure mode |
|---|---|---|---|---|---|---|---|---|---|---|---|---|---|---|---|---|
| | | | $b_w$ (mm) | $h$ (mm) | $l$ (mm) | $f_c'$ (MPa) | $f_{t,p}$ (MPa) | $G_f$ | $V_f$ (%) | $A_s$ (mm$^2$) | $f_y$ (MPa) | $f_{su}$ (MPa) | $A_p$ (mm$^2$) | $f_{py}$ (MPa) | $f_{pu}$ (MPa) | |
| (Fang et al. 2023) | B1 | T | 130 | 900 | 7450 | 140.6 | 9.2 | 16.6 | 2.0 | 402.1 | 485 | 623 | 1400.0 | 1701 | 1942 | CL |
| (El-Helou and Graybeal 2022) | B1 | T | 178 | 889 | 18890 | 173.0 | 10.4 | 16.6 | 2.0 | 0 | 0 | 0 | 4552.8 | 1690 | 1932 | CL |
| (Li et al. 2022) | B1 | T | 70 | 600 | 9800 | 128.0 | 7.1 | 20.65 | 2.0 | 905.0 | 335 | 510 | 1120 | 1701 | 1942 | GSH |
| (Shao and Billington 2022) | Ductal-vf2.0-ρ2.10 | R | 150 | 180 | 2000 | 185.8 | 10.5 | 11.2 | 2.0 | 573.0 | 454 | 778 | 0 | 0 | 0 | GSH |
| | Ductal-vf1.0-ρ2.10 | R | 150 | 180 | 2000 | 185.3 | 7.0 | 10.4 | 1.0 | 573.0 | 454 | 778 | 0 | 0 | 0 | GSH |
| | Ductal-vf2.0-ρ0.96 | R | 150 | 180 | 2000 | 185.8 | 10.5 | 11.2 | 2.0 | 253.4 | 450 | 645 | 0 | 0 | 0 | CL |
| | Ductal-vf1.0-ρ0.96 | R | 150 | 180 | 2000 | 185.3 | 7.0 | 10.4 | 1.0 | 253.4 | 450 | 645 | 0 | 0 | 0 | Transition |

Note: "T" and "R" mean bulb-tee and rectangular cross sections, respectively; "CL" means a failure after crack localization; "GSH" means a failure after gradual strain hardening of longitudinal steel reinforcement; "Transition" means a failure mode between the gradual strain hardening and crack localization.

**Appendix B.** Numerical and experimental database of flexural tests on prestressed UHPC beams.

| Name | Type | $f_{t,p}$ (MPa) | $A_s$ (mm²) | $f_{s,loc}$ (MPa) | $f_{su}$ (MPa) | $A_p$ (mm²) | $f_{con}$ (MPa) | $f_{p,loc}$ (MPa) | $f_{p,u}$ (MPa) | $d$ (mm) | $h$ (mm) | $b_w$ (mm) | $\omega$ | Ductility $\mu$ | Failure Mode |
|---|---|---|---|---|---|---|---|---|---|---|---|---|---|---|---|
| **Numerical results** | | | | | | | | | | | | | | | |
| M2@16-P10@15.2-Re | I | 9 | 402.2 | 485 | 635 | 1400 | 1395 | 1664.0 | 1942 | 820.0 | 900 | 130 | 0.50 | 1.26 | CL |
| M6@22.2-P6@15.2-Re | I | 9 | 2323.5 | 485 | 635 | 840 | 1395 | 1664.0 | 1942 | 820.0 | 900 | 130 | 0.85 | 1.71 | T |
| M10@22.7-P2@15.2-Re | I | 9 | 4044.8 | 485 | 635 | 280 | 1395 | 1664.0 | 1942 | 820.0 | 900 | 130 | 1.15 | 2.42 | GSH |
| M2@16-P10@15.2-SH50 | I | 9 | 402.2 | 485 | 535 | 1400 | 1395 | 1664.0 | 1942 | 820.0 | 900 | 130 | 0.43 | 1.18 | CL |
| M6@22.2-P6@15.2-SH50 | I | 9 | 2323.5 | 485 | 535 | 840 | 1395 | 1664.0 | 1942 | 820.0 | 900 | 130 | 0.44 | 1.67 | CL |
| M10@22.7-P2@15.2-SH50 | I | 9 | 4044.8 | 485 | 535 | 280 | 1395 | 1664.0 | 1942 | 820.0 | 900 | 130 | 0.44 | 1.77 | CL |
| M2@16-P10@15.2-SH250 | I | 9 | 402.2 | 485 | 735 | 1400 | 1395 | 1664.0 | 1942 | 820.0 | 900 | 130 | 0.57 | 1.31 | CL |
| M6@22.2-P6@15.2-SH250 | I | 9 | 2323.5 | 485 | 735 | 840 | 1395 | 1664.0 | 1942 | 820.0 | 900 | 130 | 1.26 | 2.54 | GSH |
| M10@22.7-P2@15.2-SH250 | I | 9 | 4044.8 | 485 | 735 | 280 | 1395 | 1664.0 | 1942 | 820.0 | 900 | 130 | 1.86 | 3.04 | GSH |
| M2@16-P10@15.2-fcon930 | I | 9 | 402.2 | 485 | 635 | 1400 | 930 | 1344.0 | 1942 | 820.0 | 900 | 130 | 1.29 | 1.29 | CL |
| M6@22.2-P6@15.2-fcon930 | I | 9 | 2323.5 | 485 | 635 | 840 | 930 | 1344.0 | 1942 | 820.0 | 900 | 130 | 1.32 | 2.72 | GSH |
| M10@22.7-P2@15.2-fcon930 | I | 9 | 4044.8 | 485 | 635 | 280 | 930 | 1344.0 | 1942 | 820.0 | 900 | 130 | 1.30 | 3.26 | GSH |
| M2@16-P10@15.2-ft6 | I | 6 | 402.2 | 485 | 635 | 1400 | 1395 | 1664.0 | 1942 | 820.0 | 900 | 130 | 0.76 | 1.39 | CL |
| M6@22.2-P6@15.2-ft6 | I | 6 | 2323.5 | 485 | 635 | 840 | 1395 | 1664.0 | 1942 | 820.0 | 900 | 130 | 1.28 | 3.01 | GSH |
| M10@22.7-P2@15.2-ft6 | I | 6 | 4044.8 | 485 | 635 | 280 | 1395 | 1664.0 | 1942 | 820.0 | 900 | 130 | 1.72 | 2.73 | GSH |
| M2@16-P10@15.2-ft12 | I | 12 | 402.2 | 485 | 635 | 1400 | 1395 | 1664.0 | 1942 | 820.0 | 900 | 130 | 0.38 | 1.24 | CL |
| M6@22.2-P6@15.2-ft12 | I | 12 | 2323.5 | 485 | 635 | 840 | 1395 | 1664.0 | 1942 | 820.0 | 900 | 130 | 0.64 | 1.41 | CL |
| M10@22.7-P2@15.2-ft12 | I | 12 | 4044.8 | 485 | 635 | 280 | 1395 | 1664.0 | 1942 | 820.0 | 900 | 130 | 0.86 | 1.99 | CL |
| M2@16-P10@15.2-loc0.1% | I | 9 | 402.2 | 485 | 635 | 1400 | 1395 | 1561.0 | 1942 | 820.0 | 900 | 130 | 0.76 | 1.18 | CL |
| M6@22.2-P6@15.2-loc0.1% | I | 9 | 2323.5 | 485 | 635 | 840 | 1395 | 1561.0 | 1942 | 820.0 | 900 | 130 | 1.00 | 1.36 | CL |
| M10@22.7-P2@15.2-loc0.1% | I | 9 | 4044.8 | 485 | 635 | 280 | 1395 | 1561.0 | 1942 | 820.0 | 900 | 130 | 1.20 | 2.65 | GSH |
| M2@16-P10@15.2-loc0.6% | I | 9 | 402.2 | 485 | 635 | 1400 | 1395 | 1816.0 | 1942 | 820.0 | 900 | 130 | 0.13 | 1.18 | CL |
| M6@22.2P6@15.2-loc0.6% | I | 9 | 2323.5 | 485 | 635 | 840 | 1395 | 1816.0 | 1942 | 820.0 | 900 | 130 | 0.63 | 1.29 | T |
| M10@22.7P2@15.2-loc0.6% | I | 9 | 4044.8 | 485 | 635 | 280 | 1395 | 1816.0 | 1942 | 820.0 | 900 | 130 | 1.07 | 2.49 | GSH |
| M2@16-P10@15.2-b250 | I | 9 | 402.2 | 485 | 635 | 1400 | 1395 | 1664.0 | 1942 | 820.0 | 900 | 250 | 0.26 | 1.11 | CL |
| M6@22.2P6@15.2-b250 | I | 9 | 2323.5 | 485 | 635 | 840 | 1395 | 1664.0 | 1942 | 820.0 | 900 | 250 | 0.44 | 1.33 | CL |
| M10@22.7P2@15.2-b250 | I | 9 | 4044.8 | 485 | 635 | 280 | 1395 | 1664.0 | 1942 | 820.0 | 900 | 250 | 0.60 | 1.58 | CL |

| | **Experimental results** | | | | | | | | | | | | | | | | |
|---|---|---|---|---|---|---|---|---|---|---|---|---|---|---|---|---|---|
| (Graybeal 2008) | AASHTO Type II | I | 10.3 | 0.0 | 0 | 0 | 3038.7 | 1023 | 1621.0 | 1932 | 764.0 | 914 | 152 | 0.77 | 1.00 | CL |
| (Yang et al. 2011b) | PB-65-R-S | T | 10.875 | 515.0 | 394.2 | 543.4 | 420 | 1042.2 | 1728.2 | 1942.5 | 250.0 | 300 | 50 | 1.34 | 1.46 | CL |
| | PB-00-R-S | T | 11.925 | 515.0 | 394.2 | 543.4 | 420 | 29.6 | 715.6 | 1942.5 | 250.0 | 300 | 50 | 5.62 | 3.31 | GSH |
| | PB-45-R-S | T | 11.4 | 515.0 | 394.2 | 543.4 | 420 | 775.2 | 1461.2 | 1942.5 | 250.0 | 300 | 50 | 2.49 | 2.04 | GSH |
| | FB-45-R-S-NS | T | 13.35 | 84.8 | 394.2 | 543.4 | 420 | 766.3 | 1452.3 | 1942.5 | 230.0 | 300 | 50 | 1.57 | 2.77 | CL |
| (Yang et al. 2011a) | T600NS | I | 14.5 | 0.0 | 0 | 0 | 280 | 445 | 935.0 | 1860 | 500.0 | 600 | 120 | 0.40 | 2.65 | CL |
| | T600S | I | 14.2 | 0.0 | 0 | 0 | 280 | 450 | 940.0 | 1860 | 500.0 | 600 | 120 | 0.41 | 1.46 | CL |
| | T1000S | I | 9 | 0.0 | 0 | 0 | 560 | 425 | 915.0 | 1860 | 900.0 | 1000 | 120 | 0.80 | 2.45 | CL |
| | T1300S | I | 12.1 | 0.0 | 0 | 0 | 560 | 420 | 910.0 | 1860 | 1200.0 | 1300 | 120 | 0.46 | 3.64 | CL |
| (Chen and Graybeal 2012a) | UHPC pi-girder | Pi | 9.7 | 0.0 | 0 | 0 | 2520 | 1350 | 1679.2 | 1932 | 762.0 | 840 | 170 | 0.46 | 2.00 | CL |
| (Su et al. 2020) | Box | Box | 8 | 4418.0 | 475 | 620 | 745 | 1395 | 1699.1 | 1950 | 1600.0 | 1600 | 160 | 0.67 | 3.00 | CL |
| (El-Helou and Graybeal 2022) | Bulb tee girder | T | 11 | 0.0 | 0 | 0 | 4552.8 | 1449 | 1800.0 | 1932 | 825.5 | 889 | 178 | 0.08 | 1.64 | CL |
| (Li et al. 2022) | T-girder | T | 7.1 | 905.0 | 335 | 510 | 1120 | 920 | 1410.0 | 1942 | 600.0 | 600 | 70 | 3.87 | 4.00 | GSH |
| (Fang et al. 2023) | T-girder | T | 9.2 | 402.1 | 485 | 623 | 1400 | 1395 | 1699.1 | 1942 | 820.0 | 900 | 130 | 0.40 | 1.44 | CL |
| (Xu and Deng 2014) | B1 | T | 6.5 | 56.5 | 245.2 | 381.2 | 700 | 1023 | 1513.0 | 1928.6 | 144.9 | 350 | 120 | 1.48 | 1.45 | CL |
| | B2 | T | 6.5 | 509.0 | 536 | 696.3 | 420 | 1023 | 1513.0 | 1928.6 | 217.7 | 350 | 120 | 1.41 | 3.01 | GSH |
| | B3 | T | 6.7 | 509.0 | 536 | 696.3 | 700 | 1023 | 1513.0 | 1928.6 | 211.3 | 350 | 120 | 1.92 | 3.20 | GSH |
| | B4 | T | 6.7 | 509.0 | 536 | 696.3 | 420 | 1395 | 1699.1 | 1928.6 | 217.7 | 350 | 120 | 0.85 | 3.68 | GSH |
| | B5 | T | 6.3 | 587.0 | 536 | 696.3 | 420 | 1023 | 1513.0 | 1928.6 | 225.8 | 350 | 120 | 1.54 | 3.59 | GSH |
| | B6 | T | 6.3 | 509.0 | 380.6 | 541.3 | 420 | 1023 | 1513.0 | 1928.6 | 217.7 | 350 | 120 | 1.45 | 4.86 | GSH |

Note：CL=Crack localization; GSH=Gradual strain hardening; T=Transition state.

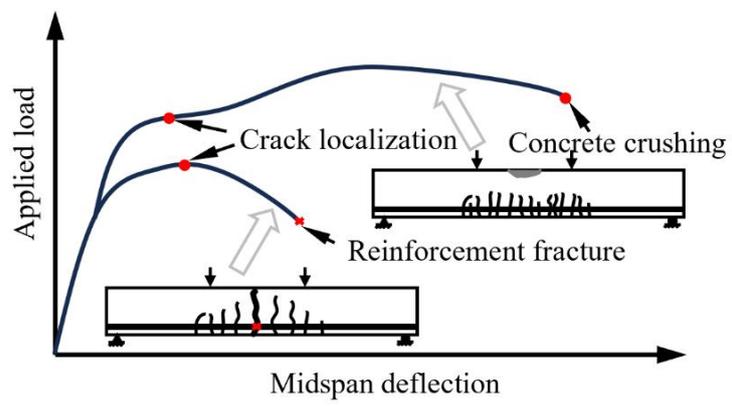

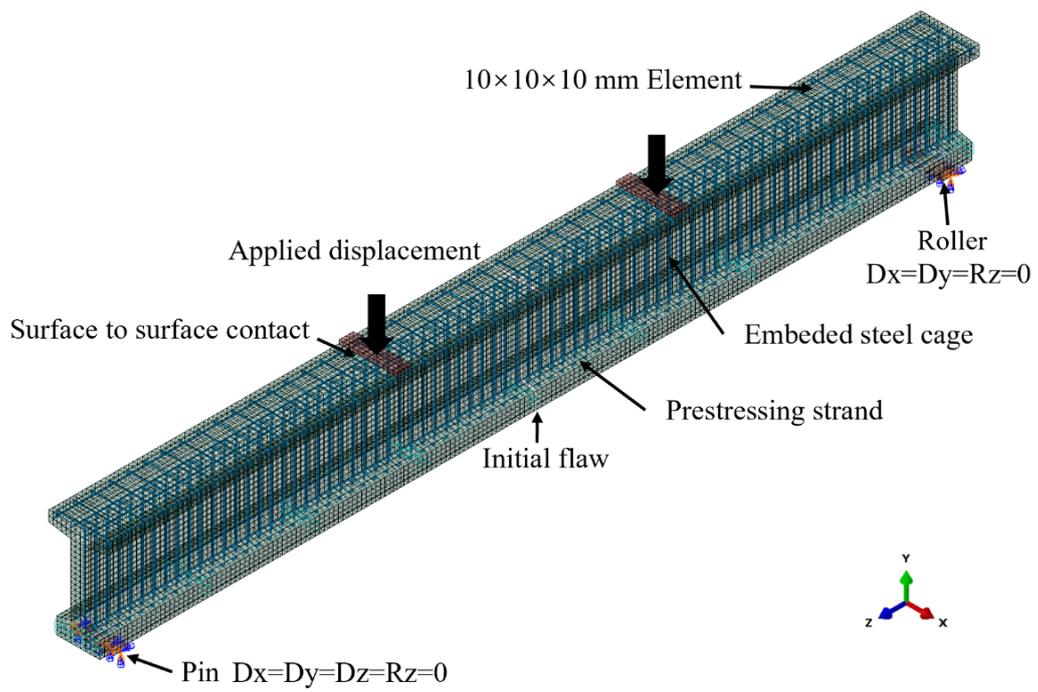

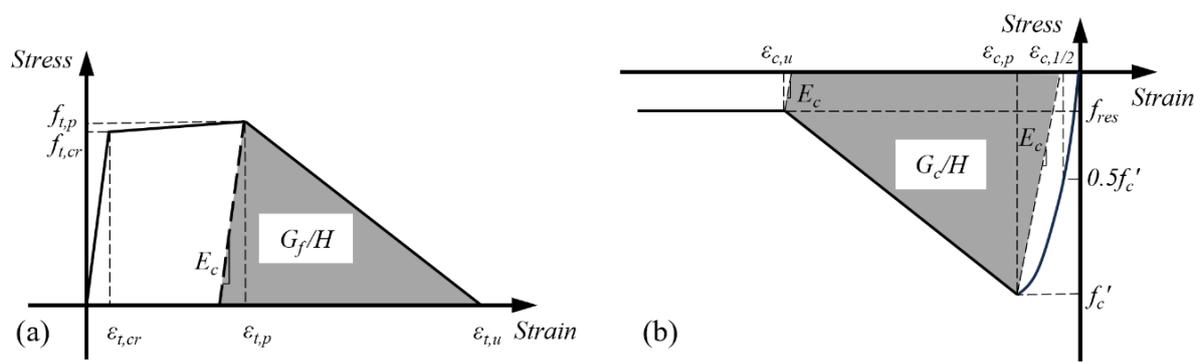

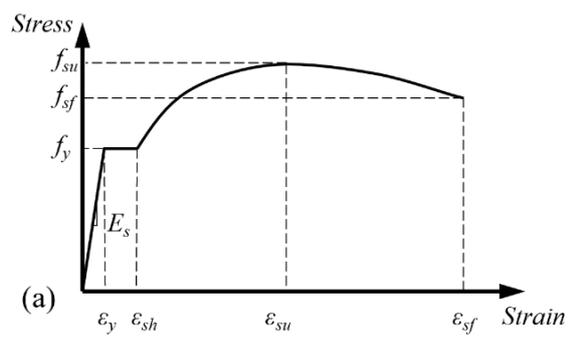 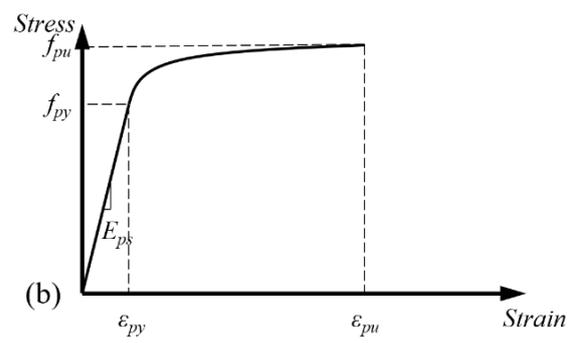

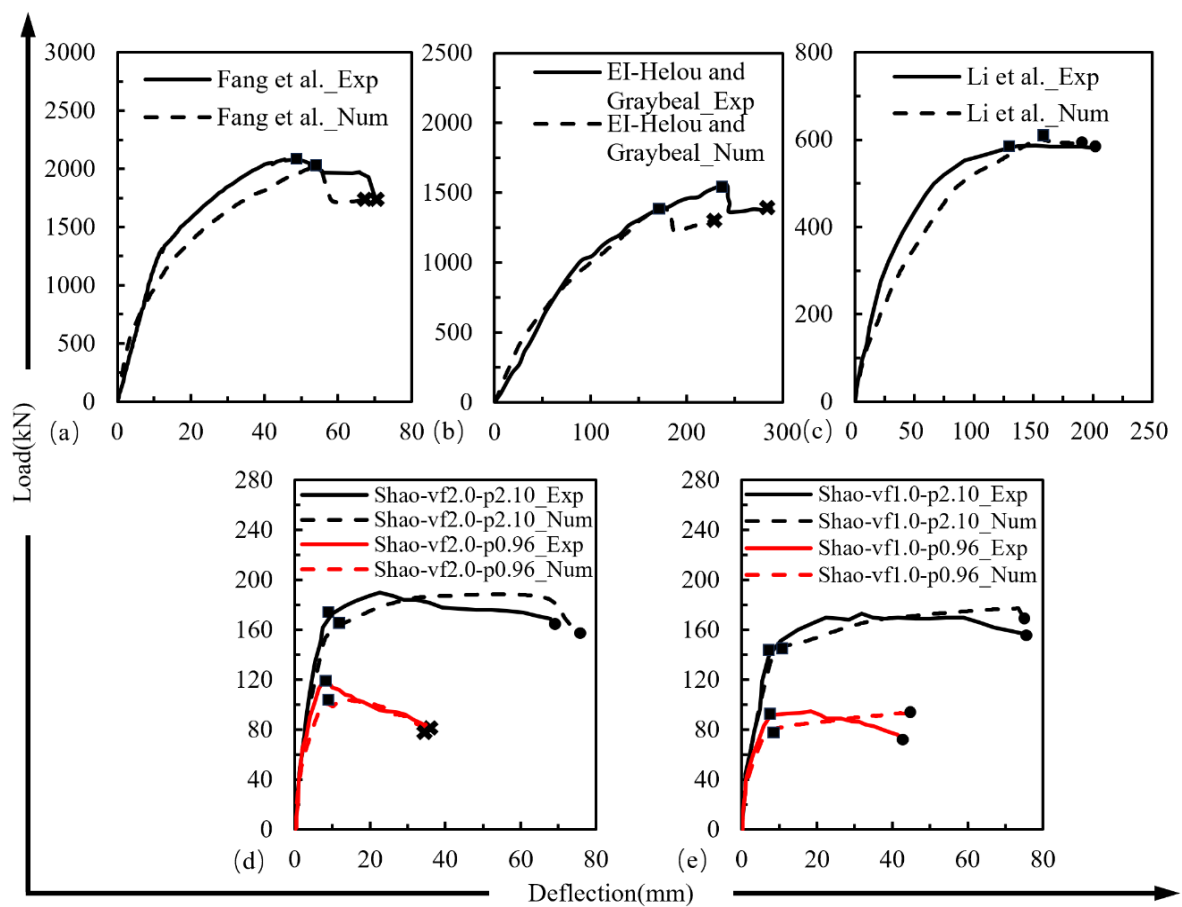

Figure6Click here to access/download;Figure;Fig 6.pdf

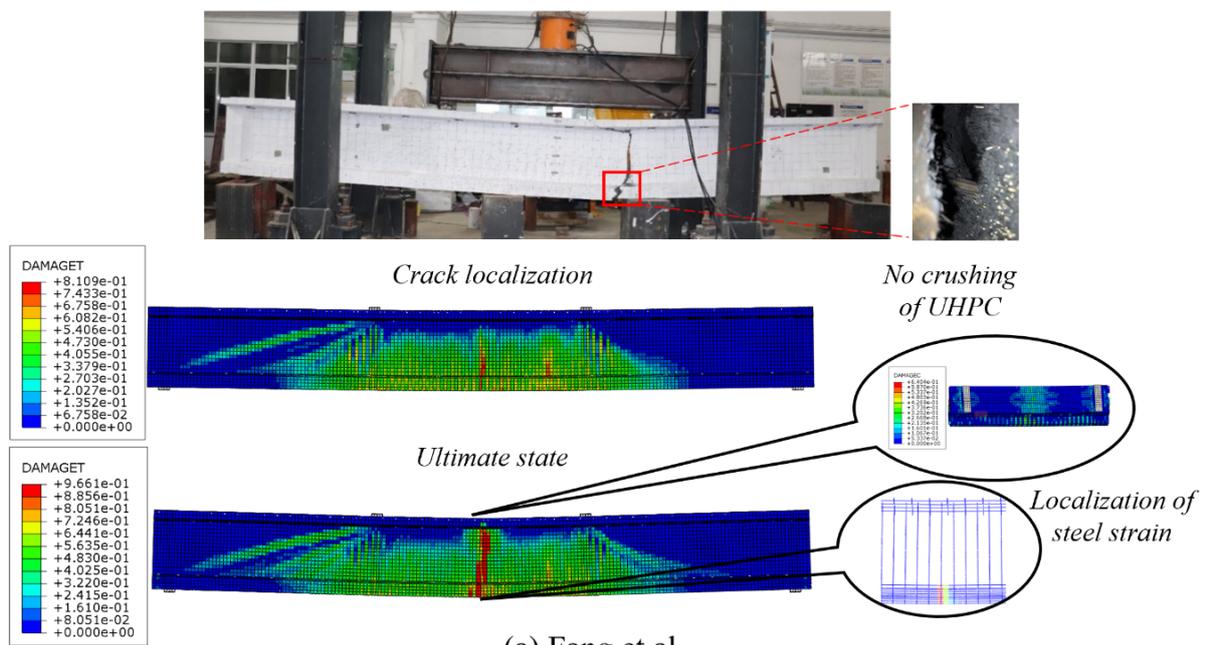

(a) Fang et al.

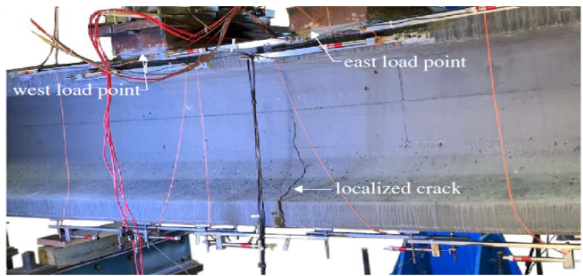

*Crack localization*

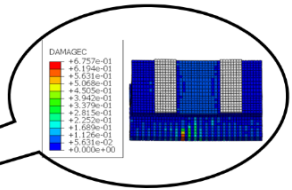

*No crushing of UHPC*

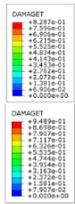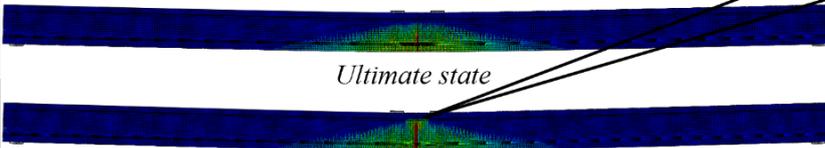

*Ultimate state*

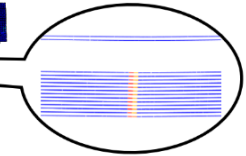

*Localization of steel strain*

(b) El-Helou and Graybeal.

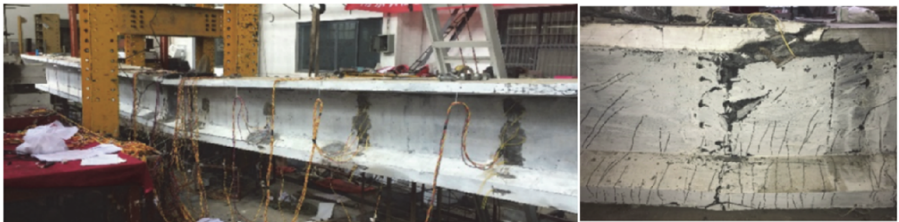
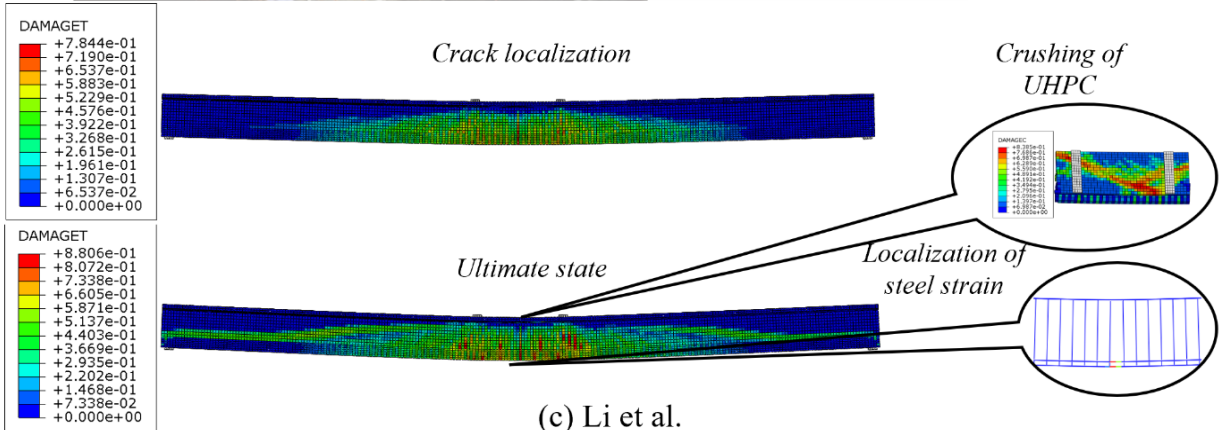

(c) Li et al.

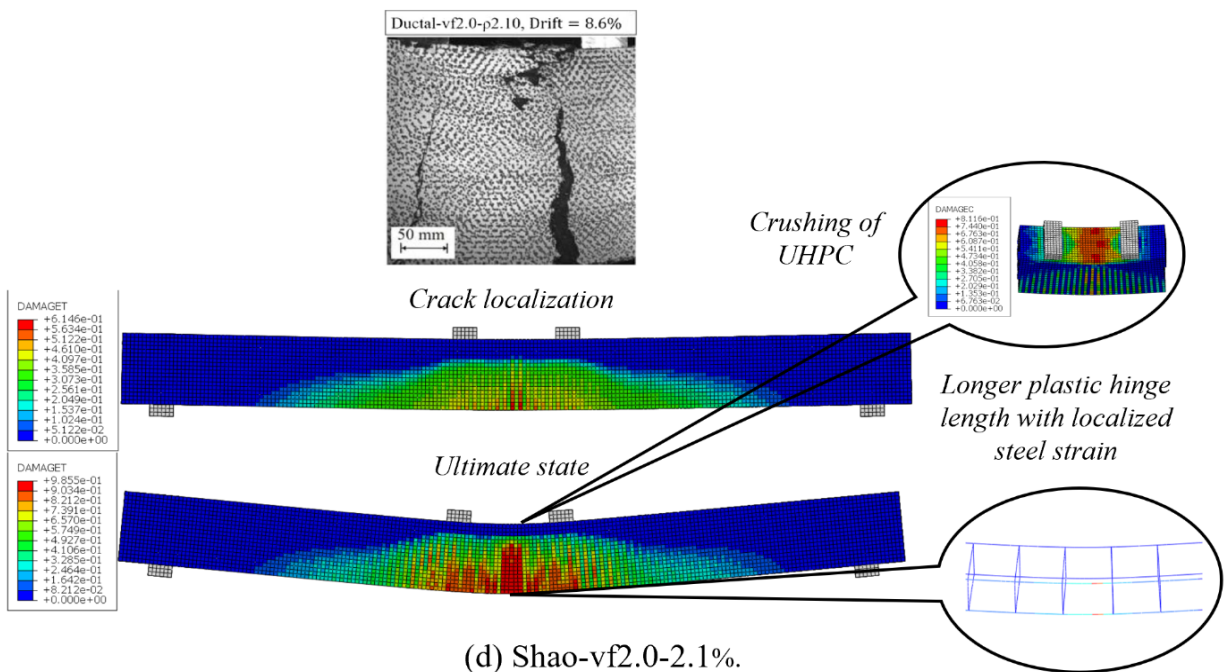

(d) Shao-vf2.0-2.1%.

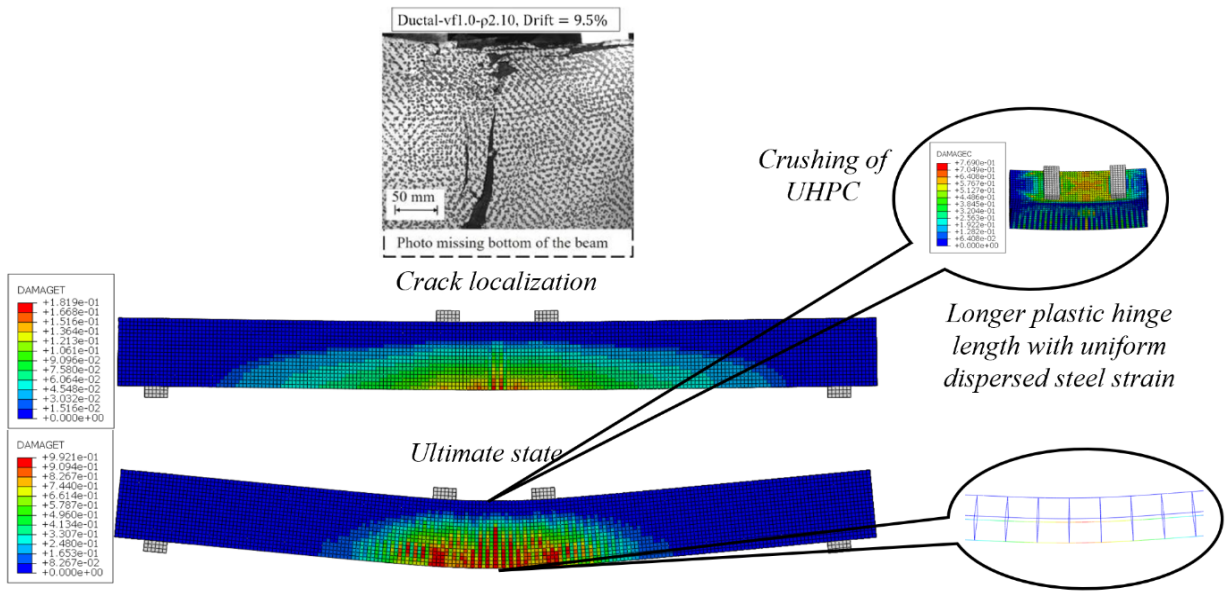

(e) Shao-vf1.0-2.1%.

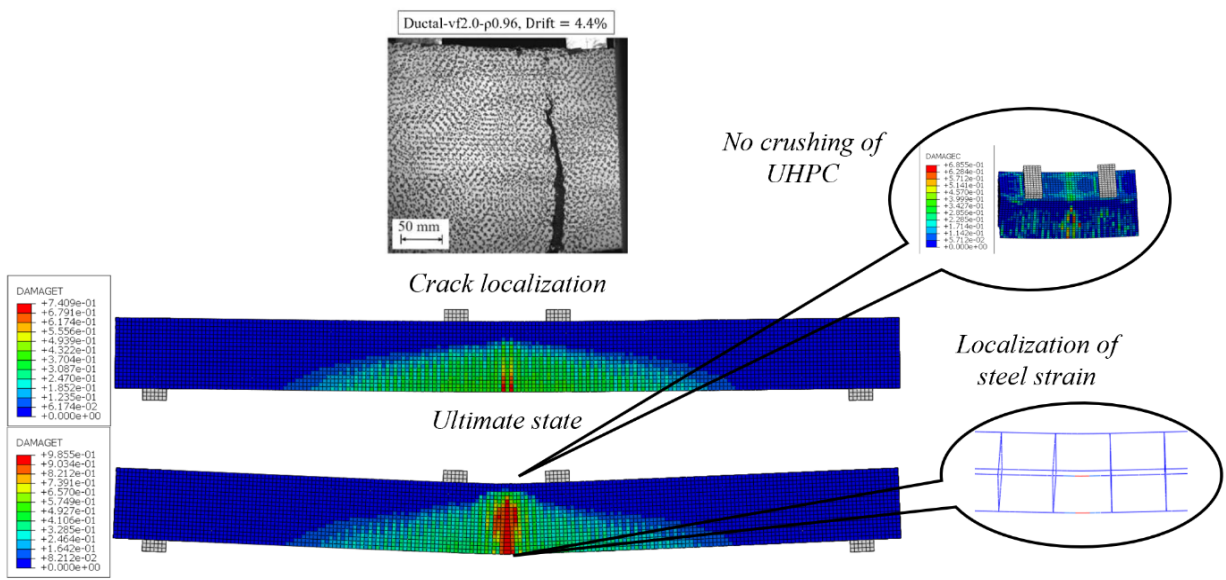

(f) Shao-vf2.0-0.96%.

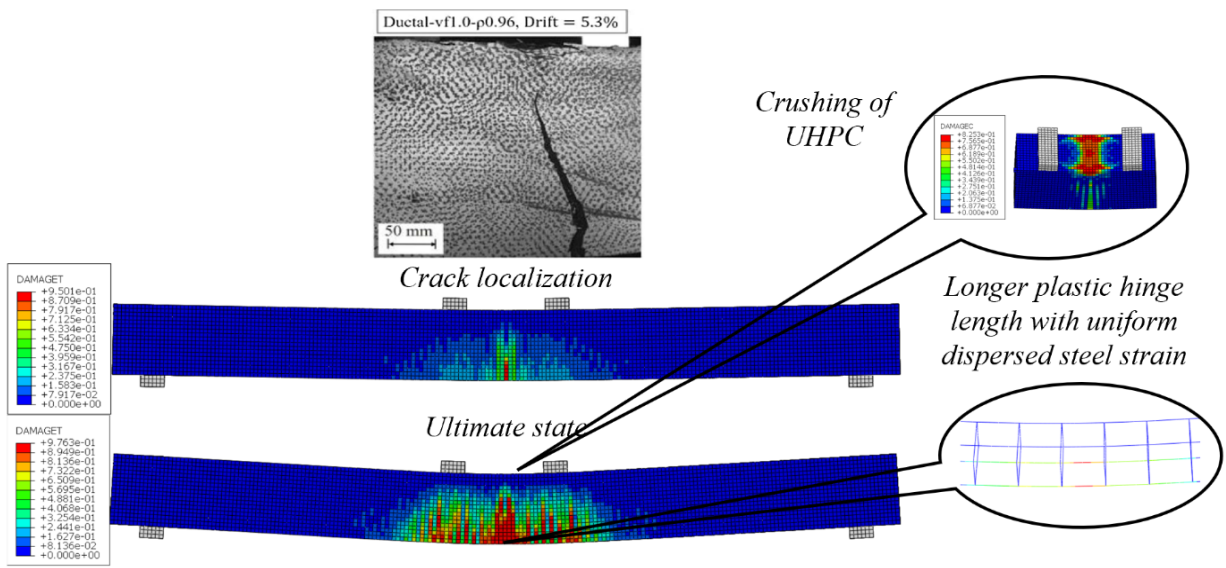

(g) Shao-vf1.0-0.96%.

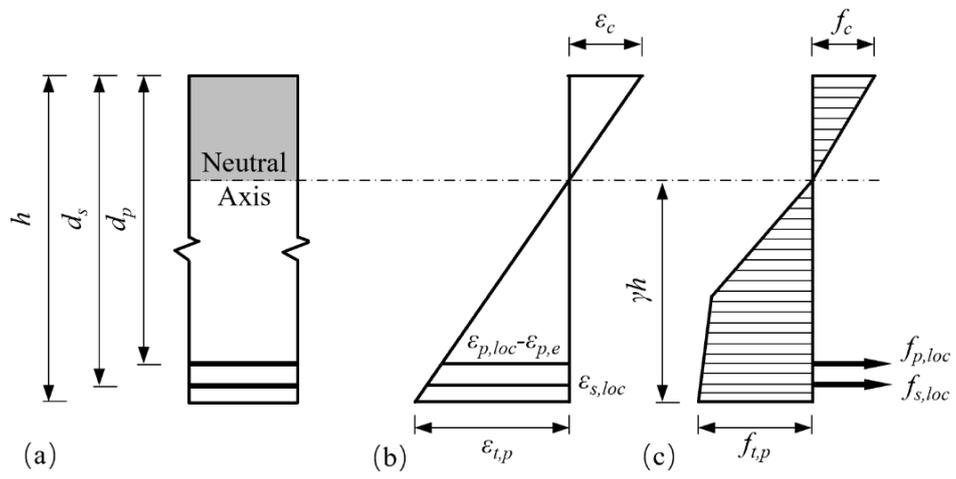

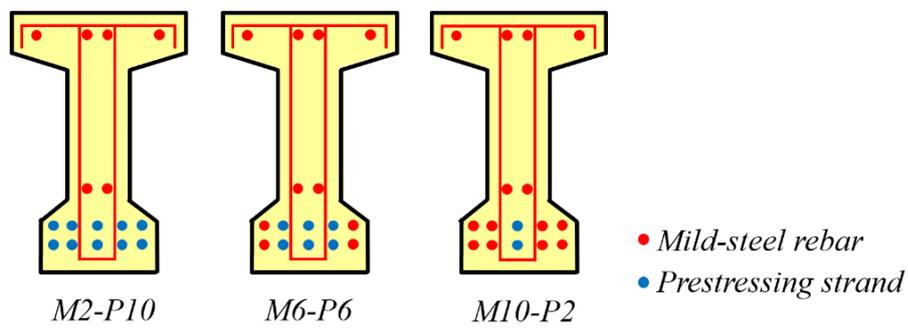

M2-P10  M6-P6  M10-P2

• *Mild-steel rebar*
• *Prestressing strand*

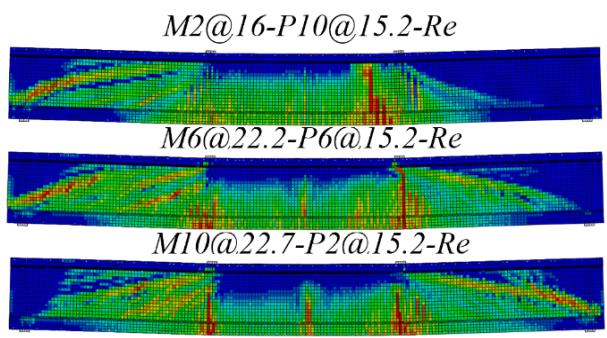

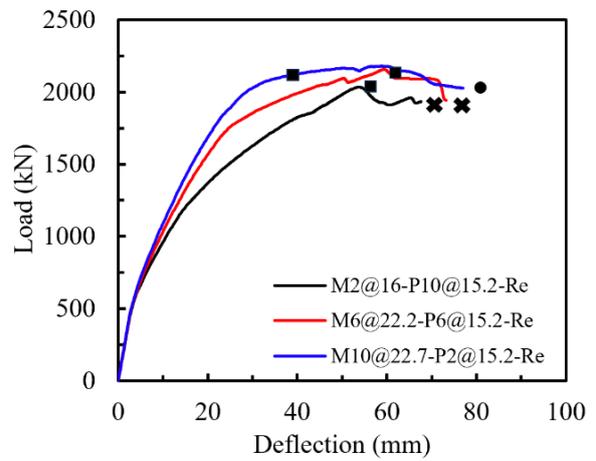

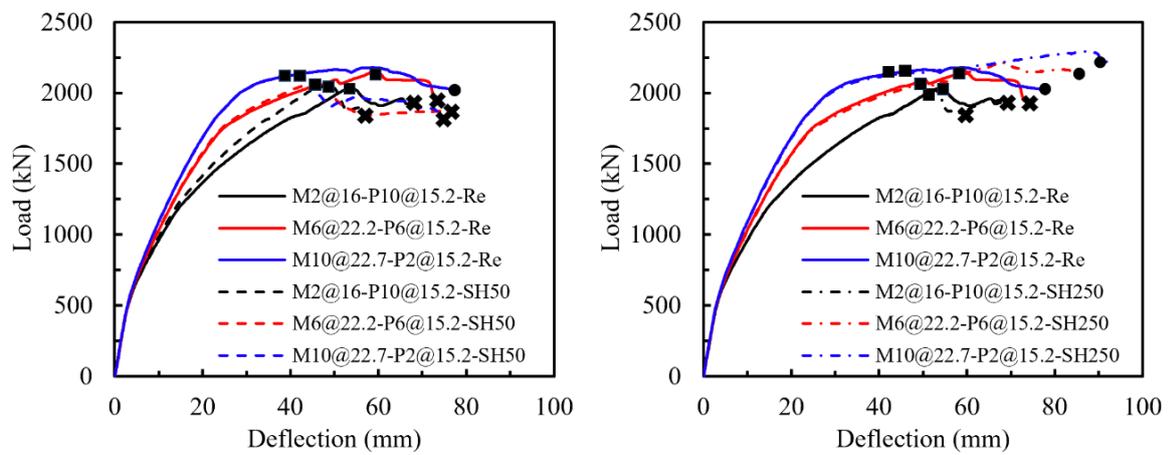

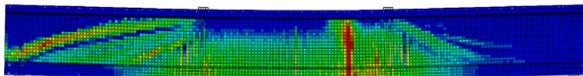
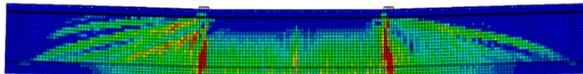
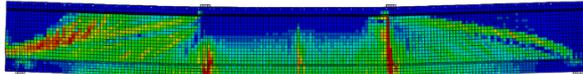
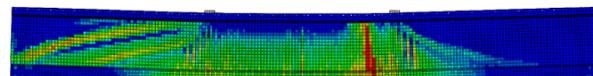
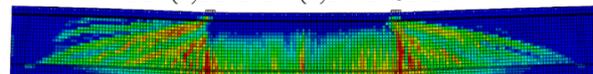
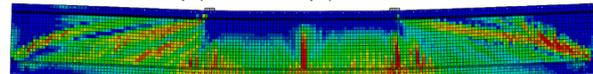

(a) $\Delta f_u$ = 50 MPa  (b) $\Delta f_u$ = 250 MPa

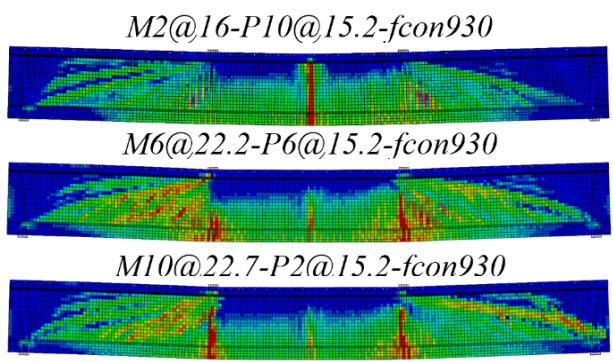
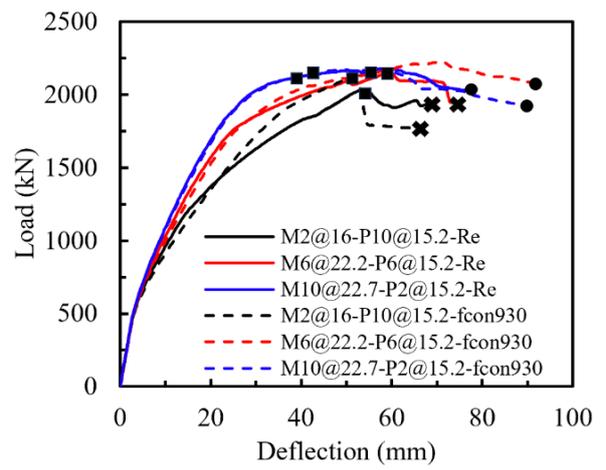

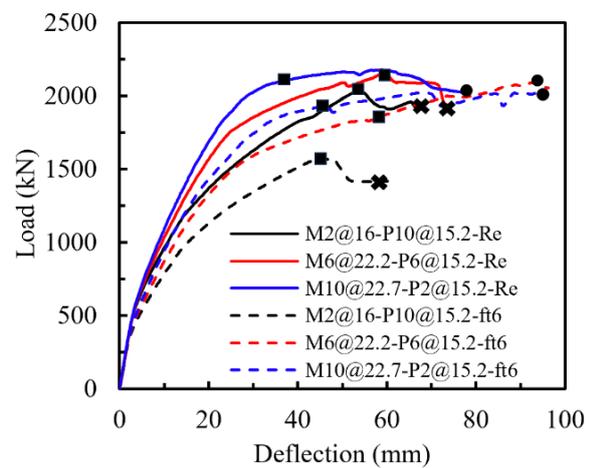
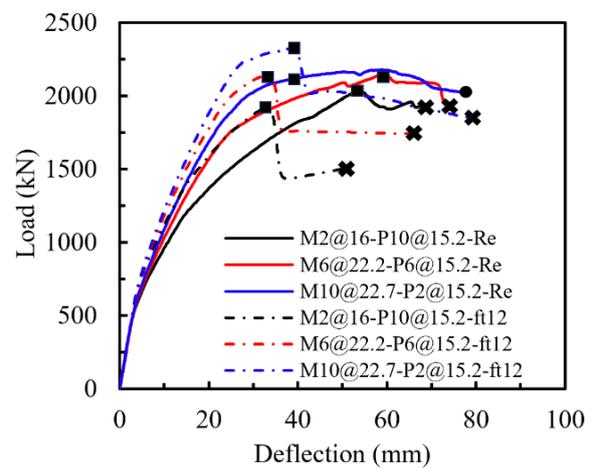
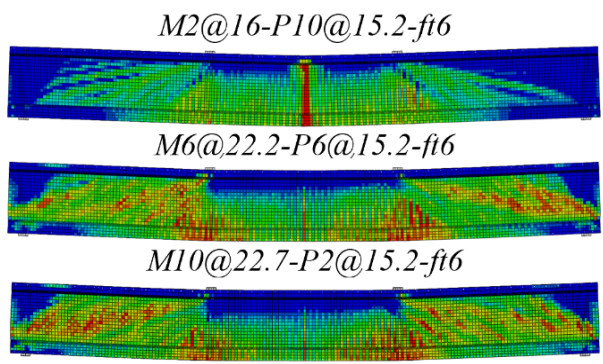
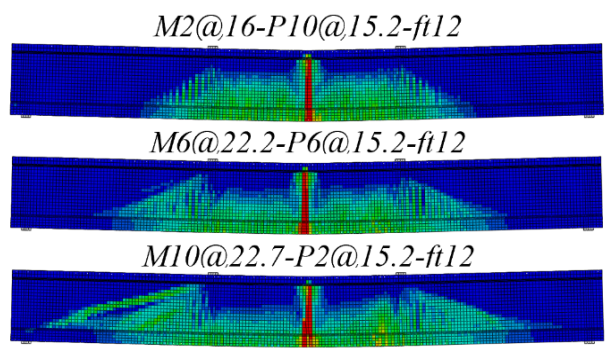

(a) $f_{t,p}$ = 6 MPa　　　　　　　　　　(b) $f_{t,p}$ = 12MPa

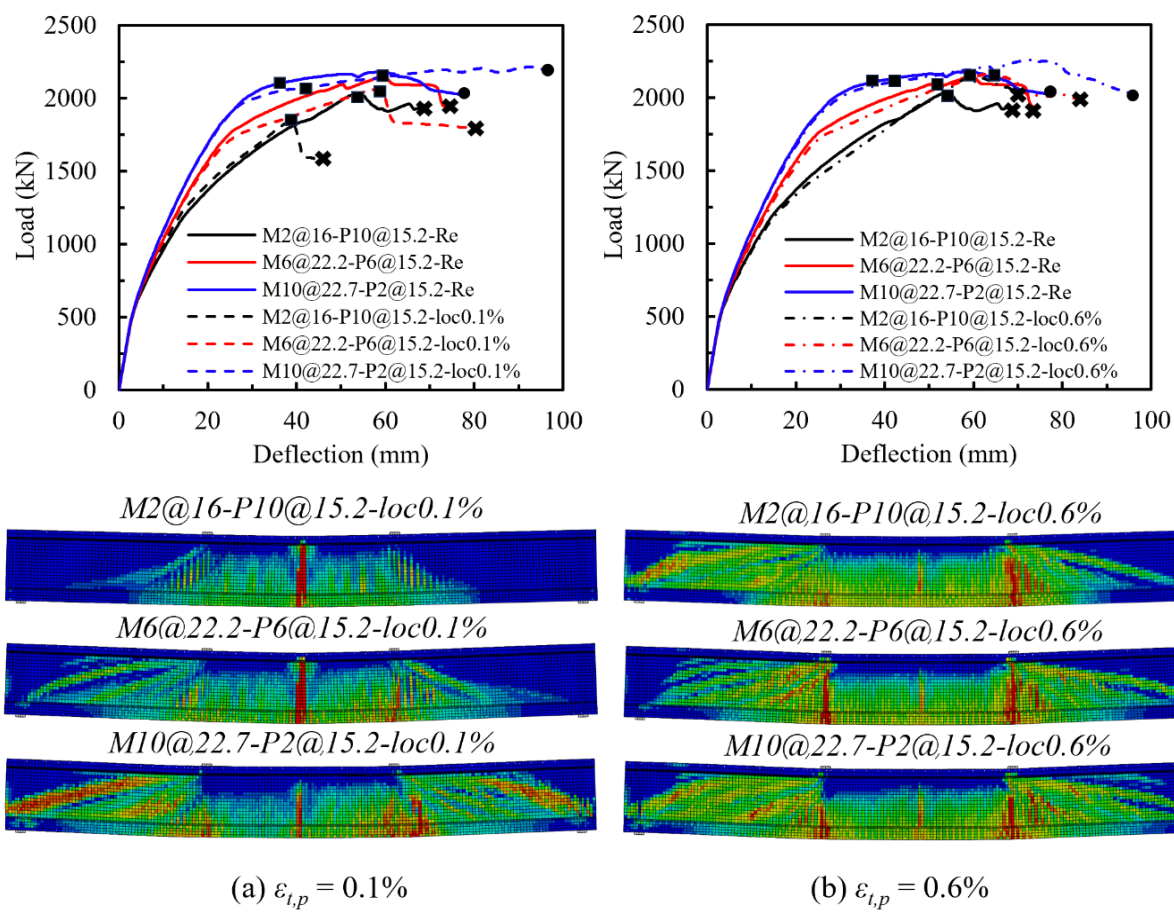

(a) $\varepsilon_{t,p} = 0.1\%$　　　　　　　　(b) $\varepsilon_{t,p} = 0.6\%$

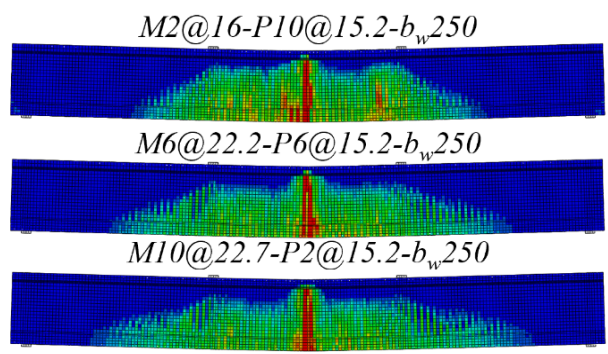
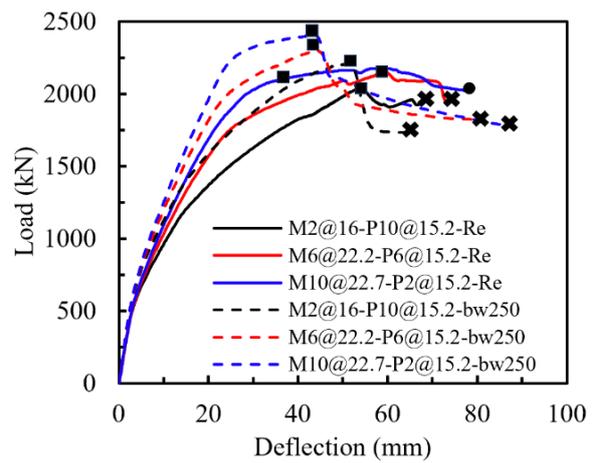

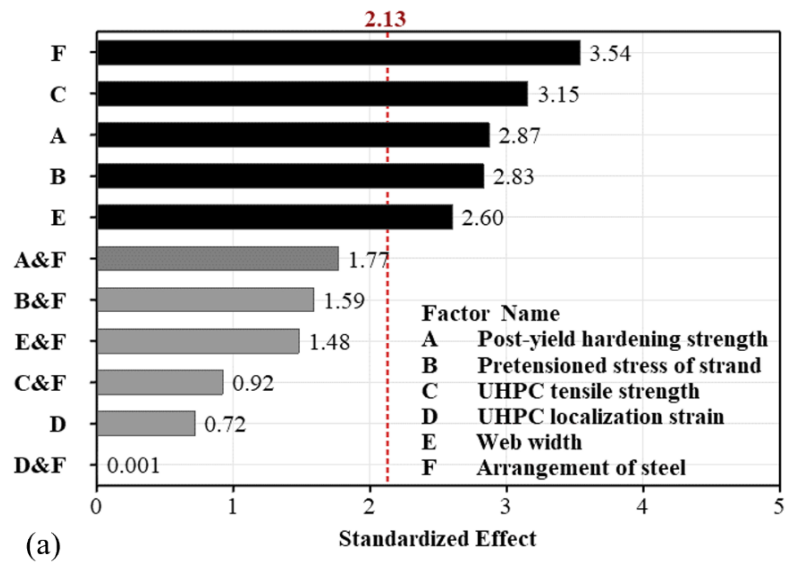
(a)

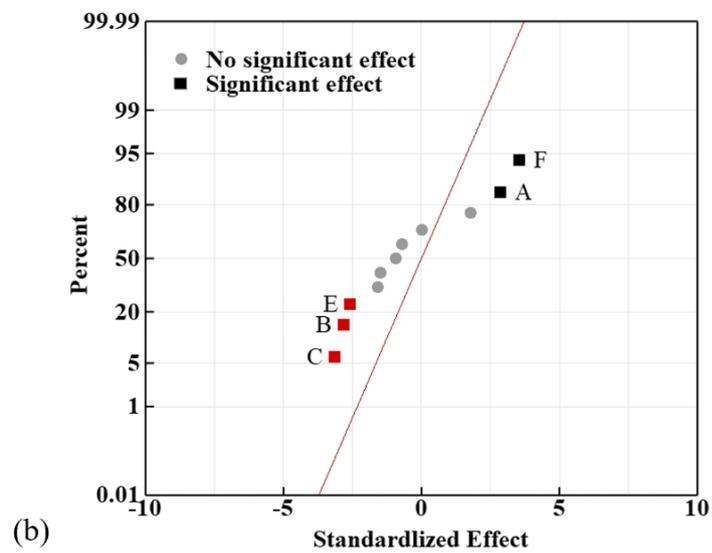
(b)

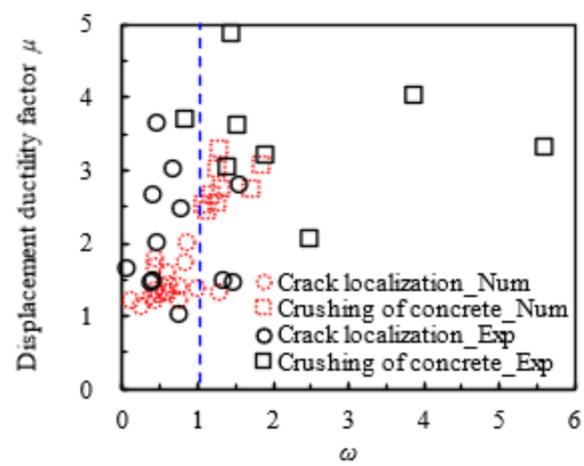

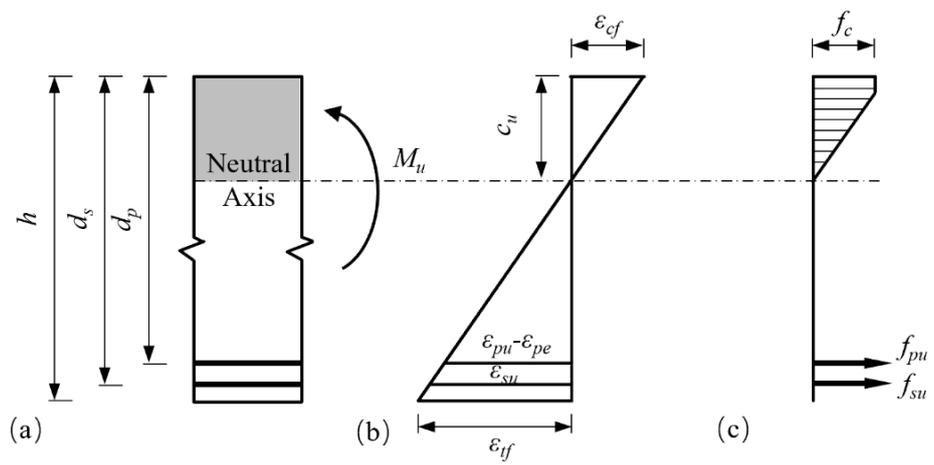

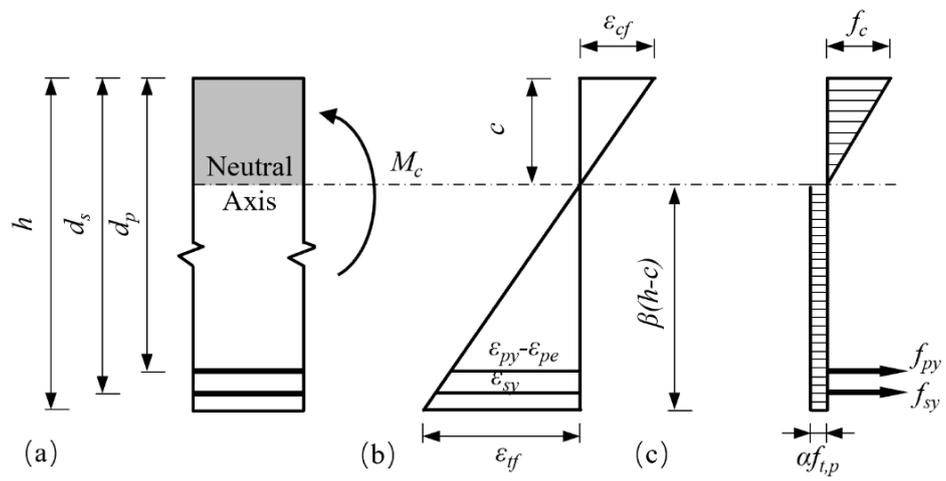

Figure Caption List

**Fig. 1.** Practical failure modes of R-UHPC beams.

**Fig. 2.** Finite-element model for typical beam.

**Fig. 3.** Idealized stress-strain model for UHPC: (a) tension; and (b) compression.

**Fig. 4.** Idealized stress-strain relationships for (a) mild steel rebar; and (b) prestressing strand.

**Fig. 5.** Comparison between the predicted and experimental load-deflection curves: (a) Fang et al.; (b) El-Helou and Graybeal.; (c) Li et al.; (d) Shao-vf2.0-2.10%. and Shao-vf2.0-0.96%; and (e) Shao-vf1.0-2.10%. and Shao-vf1.0-0.96%.; ■, ×, and ● represent the onset of crack localization, the fracture point of longitudinal steel reinforcement, and the crushing point of UHPC, respectively.

**Fig. 6.** Comparison between the predicted and experimental crack patterns: (a) Fang et al.; (b) El-Helou and Graybeal.; (c) Li et al.; (d) Shao-vf2.0-2.1%.; (e) Shao-vf1.0-2.1%.; (f) Shao-vf2.0-0.96%.; and (g) Shao-vf1.0-0.96%.

**Fig. 7.** Strain and stress distributions for prestressed UHPC beams at crack localization state: (a) beam; (b) strain; and (c) stress.

**Fig. 8.** Load-deflection responses and Tensile damage contours of specimens with different steel rebar-to-prestressing strand ratios.

**Fig. 9.** Comparison of load-deflection responses and Tensile damage contours for specimens with different post-yield hardening strength: (a) $\Delta f_u = 50$ MPa versus $\Delta f_u = 150$ MPa; and (b) $\Delta f_u = 250$ MPa versus $\Delta f_u = 150$ MPa.

**Fig. 10.** Comparison of load-deflection responses and Tensile damage contours for specimens with different pretensioned stress of strand.

**Fig. 11.** Comparison of load-deflection responses and Tensile damage contours for specimens with different tensile strengths of UHPC: (a) $f_{t,p} = 6$ MPa versus $f_{t,p} = 9$ MPa; and (b) $f_{t,p} = 12$ MPa versus $f_{t,p} = 9$ MPa.

**Fig. 12.** Comparison of Load-deflection responses and Tensile damage contours for specimens with different localization strain of UHPC: (a) $\varepsilon_{t,p} = 0.1\%$ versus $\varepsilon_{t,p} = 0.25\%$; and (b) $\varepsilon_{t,p} = 0.6\%$ versus $\varepsilon_{t,p} = 0.25\%$.

**Fig. 13.** Comparison of load-deflection responses and Tensile damage contours for specimens with different web widths.

**Fig. 14.** Influence of different factors on ductility factor based on: (a) Pareto chart of standardized effects; and (b) Normal distribution plot of standardized effects.

**Fig. 15.** Specimen displacement ductility factor $\mu$ versus $\omega$. Solid line represents experimental results, while dashed line represents numerical results.

**Fig. 16.** Strain and stress distributions for prestressed UHPC beams at ultimate state: (a) beam; (b) strain; and (c) stress.

**Fig. 17.** Strain and stress distributions for prestressed UHPC beams at crack localization state with equivalent rectangular tensile stress block: (a) beam; (b) strain; and (c) stress.